\documentclass[reqno]{amsart}
\usepackage{pictex}
\usepackage{mathrsfs}
\usepackage[colorlinks]{hyperref}
\usepackage{color}
\begin{document}

%
%

\def\labelenumi{(\theenumi)}
\newtheorem{thm}{Theorem}[section]
\newtheorem{lem}[thm]{Lemma}
\newtheorem{cor}[thm]{Corollary}
\newtheorem{add}[thm]{Addendum}
\newtheorem{prop}[thm]{Proposition}
\theoremstyle{definition}
\newtheorem{defn}[thm]{Definition}
\theoremstyle{remark}
\newtheorem{rmk}[thm]{Remark}
\newtheorem{example}[thm]{{\bf Example}}

\newcommand{\OmegaH}{\Omega/\langle H \rangle}
\newcommand{\hatOmegaHstar}{\hat \Omega/\langle H_{\ast}\rangle}
\newcommand{\SLTwoC}{\mathrm{SL}(2,\mathbf{C})}
\newcommand{\SLTwoR}{\mathrm{SL}(2,\mathbf{R})}
\newcommand{\PSLTwoC}{\mathrm{PSL}(2,\mathbf{C})}
\newcommand{\GLTwoC}{\mathrm{GL}(2,\mathbf{C})}
\newcommand{\PSLTwoR}{\mathrm{PSL}(2,\mathbf{R})}
\newcommand{\PSLTwoZ}{\mathrm{PSL}(2,\mathbf{Z})}
\newcommand{\SLTwoZ}{\mathrm{SL}(2,\mathbf{Z})}
\newcommand{\nnn}{\noindent}
\newcommand{\MCG}{{\mathcal {MCG}}}
\newcommand{\MMap}{{\bf \Phi}_{\mu}}
\def\square{\hfill${\vcenter{\vbox{\hrule height.4pt \hbox{\vrule width.4pt
height7pt \kern7pt \vrule width.4pt} \hrule height.4pt}}}$}

\newenvironment{pf}{\noindent {\it Proof.}\quad}{\square \vskip 12pt}

\title[Conditions for McShane's Identity]
{Necessary and Sufficient conditions for McShane's
Identity and variations}
\author{Ser Peow Tan, Yan Loi Wong, and Ying Zhang}
\address{Department of Mathematics \\ National University of Singapore \\
2 Science Drive 2 \\ Singapore 117543}
\email{mattansp@nus.edu.sg; matwyl@nus.edu.sg;
scip1101@nus.edu.sg}
\address{and the third author}
\address{Department of Mathematics \\Yangzhou University \\Yangzhou 225002 \\P. R. China}
\email{yingzhang@yzu.edu.cn}

\thanks{The authors are partially supported by the National University
of Singapore academic research grant R-146-000-056-112. The third
author is also partially supported by the National Key Basic
Research Fund (China) G1999075104.}

 \vskip 15pt

%
%

\begin{abstract}
Greg McShane introduced a remarkable identity for lengths of
simple closed geodesics on the once punctured torus with a
complete, finite volume hyperbolic structure. Bowditch later
generalized this and gave sufficient conditions for the identity
to hold for general type-preserving representations of a free
group on two generators $\Gamma$ to $\SLTwoC$, this was further
generalized by the authors to obtain sufficient conditions for a
generalized McShane's identity for arbitrary  (not necessarily
type-preserving) non-reducible representations. In this note we
extend the above by giving necessary and sufficient conditions for
the generalized McShane identity to hold (Akiyoshi, Miyachi and
Sakuma had proved it for type-preserving representations). We also
give a version of Bowditch's variation of McShane's identity to
once-punctured torus bundles, in the case where the monodromy is
generated by a reducible element, and an asymptotic version of
this using the generalized Markoff maps interpretation, and
provide necessary and sufficient conditions for the variations to
hold.

\end{abstract}

\maketitle

\vskip 5pt
\section{{\bf Introduction}}\label{s:intro}
\vskip 10pt

For a once punctured torus ${\mathbb T}$ with a complete, finite
volume hyperbolic structure, McShane showed
\cite{mcshane1991thesis} that the following identity holds:
\begin{eqnarray}\label{eqn:mcshanes id}
\sum_{\gamma}\frac{1}{1+e^{l(\gamma)}}=\frac{1}{2},
\end{eqnarray}
where $\gamma$ ranges over all simple closed geodesics on $\mathbb
T$, and $l(\gamma)$ is the hyperbolic length of $\gamma$ in
$\mathbb T$. This identity is independent of the hyperbolic
structure on the torus, that is, it holds for all points in the
Teichm\"uller space ${\mathcal T}_{1,1}$ of the punctured torus.
It can also be stated in terms of representations of the
fundamental group $\Gamma$ of ${\mathbb T}$ into $\SLTwoR$.

On the other hand, Bowditch gave an alternative proof of
(\ref{eqn:mcshanes id}) in \cite{bowditch1996blms} via Markoff
triples, and extended it in \cite{bowditch1998plms} to
type-preserving representations $\Gamma \rightarrow \SLTwoC$ which
satisfy certain conditions which we call here the BQ-conditions
(Definition \ref{def:BQmaps}). Akiyoshi, Miyachi and Sakuma
\cite{akiyoshi-miyachi-sakuma2004cm355} generalized this (see
Proposition 5.2 there, stated here as Theorem \ref{thm:A}) to give
(implicitly) necessary and sufficient conditions for the identity
to hold. In this note, we generalize the above to get necessary
and sufficient conditions for the generalized McShane identity to
hold for arbitrary (non-elementary) representations $\Gamma
\rightarrow \SLTwoC$ (Theorem \ref{thm:B}). Bowditch also gave
variations of McShane's identity for representations corresponding
to complete hyperbolic structures on once-punctured torus bundles
$M$ over the circle, where the monodromy of $M$ is generated by a
hyperbolic (Anosov) element of the mapping class group of the
once-punctured torus---here the mapping class group is identified
with ${\rm SL}(2, \mathbf Z)$. This was generalized in
\cite{tan-wong-zhang2004gMm} to general representations. We give
further variations in the case where the monodromy of $M$ is
generated by a reducible element (corresponding to  a parabolic
element of $\SLTwoZ$), and where some version of the BQ-conditions
holds (Theorem \ref{thm:C}). Finally, we also give an asymptotic
averaging version of this for certain generalized Markoff maps
(Theorem \ref{thm:D}).

\vskip 5pt

Let ${\mathcal X}={\rm Hom}(\Gamma, \SLTwoC)/\SLTwoC$ where the
quotient is by the conjugation action. To simplify notation, we
denote by $\rho$ instead of $[\rho]$ the conjugacy class of a
representation; there should be no confusion as we are mostly
interested in functions of the trace which are conjugacy
invariants. We denote by ${\mathcal X}_{\rm tp}\subset {\mathcal
X}$ the subset of type preserving representations; by definition,
these are the representations satisfying ${\rm tr}([\rho(a),
\rho(b)])=-2$ for any generating pair $a,b$ of $\Gamma$. We define
an equivalence relation $\sim$ on $\Gamma$ by $g \sim h$ if $g$ is
conjugate to $h$ or $h^{-1}$. Then $\Gamma/\sim$ can be identified
with the set of free homotopy classes of unoriented closed curves
on ${\mathbb T}$. Note that for $\rho \in {\mathcal X}$, the trace
function ${\rm tr}(\rho[g])$ is well-defined on the classes $[g]
\in \Gamma/\sim$.  We denote by $\hat \Omega$ the subset of
$\Gamma/\sim$ which is identified with the set of free homotopy
classes of non-trivial, non-peripheral unoriented simple closed
curves on ${\mathbb T}$.

\vskip 5pt

\nnn {\bf Conventions of Notation.}\,\, Throughout this paper, we
always assume that, for $u \in \mathbf C$,
\begin{itemize}
\item[(a)] $\sqrt{u}$ has positive real part
if $u \notin (-\infty, 0]$ and has positive imaginary part if $u
\in (-\infty, 0)$;
\item[(b)] $\log u$ ($u \neq 0$) has imaginary part in $(-\pi,
\pi]$; and
\item[(c)] $\cosh^{-1}(u)$ has real part $\ge 0$
and imaginary part $\in (-\pi, \pi]$, and has imaginary part $\in
[0, \pi]$ when real part $=0$.
\end{itemize}
The letters $\mu$ and $\tau$ are always frozen for two complex
numbers so that $\mu=\tau+2$. Implicitly, $\tau$ is for the trace
of the commutator $[A,B]:=ABA^{-1}B^{-1}$ of $A,B \in \SLTwoC$ and
$\mu$ is for ${\rm tr}^2A+{\rm tr}^2B+{\rm tr}^2AB-{\rm tr}A\,{\rm
tr}B\,{\rm tr}AB$, and the relation $\mu=\tau+2$ follows from the
well known Fricke trace identity. Consequently, the letter $\nu$
is always frozen as $\nu:=\cosh^{-1}(-\tau/2)=\cosh^{-1}(1-\mu/2)$. %

\vskip 5pt

\begin{defn}\label{def:hfunction}
Define the function $h: {\mathbf C} \backslash \{0\} \rightarrow
{\mathbf C}$ by
\begin{eqnarray}\label{eqn:h(x)=}
h(x)=\frac{1}{2}\bigg(1-\sqrt{1-\frac{4}{x^2}}\,\bigg).
\end{eqnarray}
\end{defn}

Then Bowditch's extension and reformulation of (\ref{eqn:mcshanes
id}), stated in terms of representations, is as follows.

\vskip 5pt

\begin{thm}\label{thm:BM} {\rm (Theorem 3 \cite{bowditch1998plms})}
Suppose $\rho \in {\mathcal X}_{\rm tp}$ satisfies

{\rm(i)}\,\, ${\rm tr}\rho[g] \not\in [-2,2]$ for all $[g] \in
\hat \Omega$; and

{\rm(ii)} $|{\rm tr}\rho[g]| \le 2$ for only finitely many classes
$[g] \in \hat \Omega$.

\nnn Then
\begin{eqnarray}\label{eqn:Bowditch}
\sum_{[g] \in \hat \Omega}h\big( {\rm tr}\rho[g]\big)=\frac{1}{2},
\end{eqnarray}
where the sum converges absolutely.
\end{thm}

We call conditions (i) and (ii) in Theorem \ref{thm:BM} the
BQ-conditions. Note that the summands in (\ref{eqn:Bowditch}) are
equal to the summands in (\ref{eqn:mcshanes id}) for
representations coming from complete finite area hyperbolic
structures on ${\mathbb T}$. Bowditch has conjectured that the
representations satisfying the BQ-conditions are precisely the
quasifuchsian representations (Conjecture A in
\cite{bowditch1998plms}). Note also that because $h(x)$ is an even
function, the theorem holds for representations into ${\rm PSL}(2,
\mathbf C)$ as well (the trace of the commutator is well-defined).

Condition (i) implies that there are no elliptics or accidental
parabolics in $\rho(\hat \Omega)$ and condition (ii) is clearly
necessary for the sum in (\ref{eqn:Bowditch}) to converge. It
turns out these conditions are almost necessary; it is easy to see
that it is necessary that there are no elliptics (otherwise there
are infinitely many $[g]\in \hat \Omega$ with  $|{\rm
tr}(\rho[g])|<K$ for some $K>0$ and the sum diverges), but it
turns out that one can allow for (a finite number of) accidental
parabolics and (\ref{eqn:Bowditch}) will still hold. We have the
following result of Akiyoshi, Miyachi and Sakuma
\cite{akiyoshi-miyachi-sakuma2004cm355}.

\begin{thm}\label{thm:A}{\rm
(\cite{akiyoshi-miyachi-sakuma2004cm355}, Proposition 5.2)}\,
Suppose $\rho \in {\mathcal X}_{\rm tp}$. Then
\begin{eqnarray*}
\sum_{[g] \in \hat \Omega}h\big({\rm tr}\rho[g]\big)=\frac{1}{2},
\end{eqnarray*}
and the sum converges absolutely, if and only if

{\rm (i$'$)}\, ${\rm tr}\rho[g] \not\in (-2,2)$ for all $[g] \in
\hat \Omega$, and

{\rm (ii)} $|{\rm tr}\rho[g]| \le 2$ for only finitely many
classes $[g] \in \hat \Omega$.
\end{thm}

We will see later (Example \ref{ex:AA}) that conjecturally, the
extra representations included in Theorem \ref{thm:A} come from
the points in the Maskit embedding of ${\mathcal T}_{1,1}$ (see
\cite{keen-series1993t} for background and a discussion of the
geometry of this embedding), and the cusp points on the boundary
of the embedding. In particular, there are one or two accidental
parabolics in these extremal cases.

\vskip 10pt

We turn now to more general representations. A representation
$\rho \in {\mathcal X}$ is said to be a $\tau$-{\bf
representation}, where $\tau \in \mathbf C$, if for some (hence
every) pair of free generators $a, b \in \Gamma$, ${\rm
tr}\rho([a,b]) = \tau$. Denote by ${\mathcal X}_{\tau} \subset
{\mathcal X}$ the set of $\tau$-representations (so ${\mathcal
X}_{\rm tp}={\mathcal X}_{-2}$). Note that this makes sense for
representations into ${\rm PSL}(2, \mathbf C)$ as well since the
trace of the commutator is well-defined. It is well known that
$\rho(\Gamma)$ is elementary if and only if $\tau=2$, which
corresponds to reducible representations, and a different type of
geometry; this will be explored in a future paper. Hence, for now,
we will only be interested in $\tau$-representations with $\tau
\neq 2$. To state the generalized version of McShane's identity
formulated in terms of representations, we need the following
definition.

\begin{defn}\label{defn:frak h}
For a fixed $\tau \neq \pm 2$, let $\nu=\cosh^{-1}(-\tau/2)$.
Define a function ${\mathfrak h}_{\tau}: {\mathbf C} \backslash
\{\pm \sqrt{\tau+2}\} \rightarrow {\mathbf C}$ by
\begin{eqnarray}\label{eqn:frak h(x)=}
{\mathfrak h}_{\tau}(x)=\log
\frac{1+(e^{\nu}-1)h(x)}{1+(e^{-\nu}-1)h(x)},
\end{eqnarray}
where $h(x)$ is given as in (\ref{eqn:h(x)=}).
\end{defn}

\vskip 5pt

\begin{rmk}
It can be verified by an elementary calculation that
\begin{eqnarray*}
{\mathfrak h}_{\tau}\big( {\rm tr}\rho[g]\big)=2
\tanh^{-1}\left(\frac{\sinh\nu}{\cosh\nu+e^{l(\rho[g])}} \right),
\end{eqnarray*}
where $l(\rho[g])$ is the complex translation length of $\rho[g]$.
This was the version used in
\cite{tan-wong-zhang2004cone-surfaces}, and has a more geometric
flavor.
\end{rmk}

\vskip 5pt

\begin{thm}\label{thm:B}
Let $\rho: \Gamma \rightarrow {\rm SL}(2, \mathbf C)$ be a
$\tau$-representation, where $\tau \neq \pm 2$, and let
$\nu=\cosh^{-1}(-\tau/2)$. Then
\begin{eqnarray} \label{eqn:TWZ}
\sum_{[g] \in \hat \Omega}{\mathfrak h}_{\tau} \big( {\rm
tr}\rho[g]\big)=\nu \mod 2\pi i,
\end{eqnarray}
and the sum converges absolutely, if and only if

{\rm (i$'$)}\, ${\rm tr}\rho[g] \not\in (-2,2)$ for all $[g] \in
\hat \Omega$, and

{\rm (ii)}  $|{\rm tr}\rho[g]| \le 2$ for only finitely many
classes $[g] \in \hat \Omega$.
\end{thm}

\vskip 5pt

\begin{rmk} Note that ${\mathfrak h}_{\tau}$ is not
defined at $\pm \sqrt{\tau+2}$, but this does not cause any
essential restriction since it was shown in
\cite{tan-wong-zhang2004gMm} that if ${\rm tr}\rho[g]=\pm
\sqrt{\tau+2}$ for some $[g] \in \hat\Omega$, then part(ii) of the
BQ-conditions is not satisfied. It was also shown in
\cite{tan-wong-zhang2004gMm} that the BQ-conditions are sufficient
for (\ref{eqn:TWZ}) to hold with absolute convergence; Theorem
\ref{thm:B} is an extension of the result there. To get a feeling
for the theorem, the reader should consider representations
arising from hyperbolic structures on the three-holed sphere (pair
of pants) with geodesic boundary---they satisfy the
BQ-conditions---and then consider deformations where some, or all
of the boundary components degenerate to cusps; see Example
\ref{ex:BB} (these satisfy conditions (i$'$) and (ii) in Theorem
\ref{thm:B}). See also \cite{goldman2003gt} for a more detailed
analysis of all such real representations.
\end{rmk}

To state our next result, we consider the action of the mapping
class group of $\mathbb T$, ${\mathcal {MCG}} \cong \SLTwoZ$, and
its induced action  on $\Gamma$ and ${\mathcal X}_{\tau}$. Recall
that $\MCG$ is the group of isotopy classes of diffeomorphisms of
${\mathbb T}$ which fixes the puncture (that is, each element of
$\MCG$ fixes a neighborhood of the puncture pointwise). Every $H
\in \MCG$ induces an automorphism $H_{\ast}$ of $\Gamma$ and $H$
acts on ${\mathcal X}_{\tau}$ by
$$H(\rho)(g)=\rho(H_{\ast}(g))$$
for $\rho \in {\mathcal X}_{\tau}$ and $g \in \Gamma$. Bowditch
\cite{bowditch1997t} studied representations $\rho \in {\mathcal
X}_{\rm tp}$ stabilized by a cyclic subgroup $\langle H \rangle
<\MCG \cong \SLTwoZ$ generated by a hyperbolic element and proved
a variation of the McShane's identity. This was subsequently
generalized in \cite{tan-wong-zhang2004gMm}. The result (again
restated here in terms of representations) is as follows.

\vskip 5pt

\begin{thm}{\rm (Theorem A, \cite{bowditch1997t}; see also
\cite{tan-wong-zhang2004gMm})}\label{thm:BowditchT} Suppose that
$\rho \in {\mathcal X}_{\tau}$ {\rm(}$\tau \neq 2${\rm)} is
stabilized by $\langle H \rangle <\MCG \cong \SLTwoZ$, where $H$
is a hyperbolic element of $\SLTwoZ$. Suppose

{\rm (i)}\,\,  ${\rm tr}\rho[[g]] \not\in [-2,2]$ for all $[[g]]
\in {\hatOmegaHstar}$, and

{\rm (ii)}   $|{\rm tr}\rho[[g]]| \le 2$ for only finitely
many classes $[[g]] \in {\hatOmegaHstar}$. \\
Then
\begin{eqnarray}\label{eqn:hyperbolicstab}
\sum_{[[g]] \in \hatOmegaHstar}{\mathfrak h}\big( {\rm
tr}\rho[[g]]\big)=0 \mod 2\pi i,
\end{eqnarray}
where if $\tau=-2$ then ${\mathfrak h}(x)=h(x)$, and if $\tau \neq
-2$ then ${\mathfrak h}(x)={\mathfrak h}_{\tau}(x)$. Moreover, the
sum converges absolutely.
\end{thm}

Using essentially the same arguments, we can relax the condition
(i) to:

\vskip 3pt

{\rm (i$'$)}  ${\rm tr}\rho[[g]] \not\in (-2,2)$ for all $[[g]]
\in {\hatOmegaHstar}$, \vskip 3pt

\nnn as in Theorems \ref{thm:A} and \ref{thm:B}. However, it is
not difficult to see that for fixed $\tau$, the set of $\rho \in
{\mathcal X}_\tau$ stabilized by $\langle H\rangle$ is finite (so
there is a rigidity-type property), so it is possible that this
does not give any new examples of representations satisfying this
variation of the identity. Instead, we prove the following version
of the above where the stabilizer $\langle H \rangle$ is generated
by a reducible (parabolic) element $H$, and the deformation space
is more interesting. In this case the induced action $H_{\ast}$ of
$H$ on $\Gamma$ is given by $H_{\ast}(a)=a^{\pm1}$,
$H_{\ast}(b)=a^{ N}b^{\pm 1}$, $N \in {\mathbf Z}$, for some fixed
generating pair $a,b$ of $\Gamma$.

\begin{thm}\label{thm:C}
Suppose that $\rho \in {\mathcal X}_{\tau}$ {\rm(}$\tau \neq
2${\rm)} is stabilized by $\langle H \rangle <\MCG \cong \SLTwoZ$,
where $H$ is a parabolic element of $\SLTwoZ$. Let $[g_0]$ be the
unique element in $\hat \Omega$ fixed by $H_{\ast}$. Then
\begin{eqnarray}\label{eqn:parabolicstab}
\sum_{[[g]] \in \hatOmegaHstar-[[g_0]]}{\mathfrak h}\big( {\rm
tr}\rho[[g]]\big)=0 \mod 2\pi i,
\end{eqnarray}
where ${\mathfrak h}=h$ if $\tau=-2$, and ${\mathfrak
h}={\mathfrak h}_{\tau}$ if $\tau \neq -2$, and the sum converges
absolutely, if and only if

{\rm (i)}\,\, ${\rm tr}\rho[[g]] \not\in (-2,2)$ for all $[[g]]
\in \hatOmegaHstar-[[g_0]],$ and

{\rm (ii)} $|{\rm tr}\rho[[g]]| \le 2$ for only finitely many
classes $[[g]] \in {\hatOmegaHstar}$.
\end{thm}

\vskip 5pt

%

The above can also be stated in terms of once punctured torus
bundles over the circle with monodromy a reducible element of
$\MCG$. In this case $\hatOmegaHstar$ can be identified with the
set of free homotopy classes of non-trivial, non-peripheral simple
closed curves on the fibre, and $\hatOmegaHstar-[[g_0]]$ with the
subset which excludes the class of the fixed curve. Note that
there is no generalization of Theorem B of \cite{bowditch1997t} in
this case, because $\hatOmegaHstar-[[g_0]]$ is essentially
one-sided. It is also not difficult to see that, for ${\mathcal
X}_{\rm tp}$, the remaining case of representations stabilized by
a subgroup generated by an elliptic element $H$ gives rise to
either representations which satisfy the BQ-conditions, or the
representation corresponding to the trivial Markoff map, that is,
where ${\rm tr}\rho(g)=0$ for all $g \in \Gamma$. Hence Theorems
\ref{thm:BowditchT} and \ref{thm:C} together with the observation
above give a rather complete answer to the cases of
representations in ${\mathcal X}_{\rm tp}$ stabilized by either
finite or (virtually) cyclic subgroups. In a future paper
\cite{tan-wong-zhang2004endinvariants} we will explore
representations in ${\mathcal X}_{\rm tp}$ stabilized by subgroups
which are not finite or virtually cyclic, and also the end
invariants of a generalized Markoff map $\phi$. It turns out that
the extended Bowditch's conditions are more useful there. We defer
the statement of our final result (Theorem \ref{thm:D}) to the
next section as it is more conveniently stated in terms of
generalized Markoff maps.

\vskip 5pt

We give a brief discussion of the proofs. Following Bowditch
\cite{bowditch1998plms} \cite{bowditch1997t}, the results above
are reformulated (see \S \ref{sec:Markoffmaps}) in terms of
Markoff maps (for Theorem \ref{thm:A}); generalized Markoff maps
(for Theorem \ref{thm:B}); and generalized Markoff maps stabilized
by a parabolic element of $\PSLTwoZ$ (for Theorem \ref{thm:C}).
The basic idea in \cite{bowditch1998plms}, extended in
\cite{tan-wong-zhang2004gMm}, was that the BQ-conditions imply
that the (generalized) Markoff maps have Fibonacci growth either
on $\Omega$ (for Theorems \ref{thm:A} and \ref{thm:B}) or on
suitable branches of the tree $\Sigma$ (for Theorem \ref{thm:C}),
which was sufficient to give the absolute convergence of the sums
in Theorems \ref{thm:A}, \ref{thm:B} and \ref{thm:C}. The actual
sum was then computed by a tricky computation, using the
$\psi$-functions $\frac{x}{yz}$, $\frac{y}{zx}$ and $\frac{z}{xy}$
associated to a Markoff triple $(x,y,z)$, or the functions
$\Psi(y,z,x),\Psi(z,x,y),\Psi(x,y,z)$ associated to a generalized
Markoff triple $(x,y,z)$. The existence of accidental parabolics
(regions where the generalized Markoff map takes values $\pm 2$)
means that the associated generalized Markoff map does not have
Fibonacci growth. However, it turns out that the condition of
Fibonacci growth is stronger than is necessary for the absolute
convergence of the sums and the computation of their values. By
analyzing the behavior of a generalized Markoff map on the regions
surrounding a region with value $\pm 2$, we show that the absolute
convergence of the sums with the resulting values in Theorems
\ref{thm:A}, \ref{thm:B} and \ref{thm:C} still holds in the
presence of a finite number of accidental parabolics, so that the
BQ-conditions can be relaxed (see \S \ref{sec:proofs}, and also
\cite{akiyoshi-miyachi-sakuma2004cm355} where the proof follows
along very similar lines). There are some subtleties involved in
the proof of Theorem \ref{thm:C} for some special cases; these are
taken care of in Lemma \ref{lem:x^2=mu}. The proof of Theorem
\ref{thm:D} is similar to that for Theorem \ref{thm:C}. Because of
the non-periodicity, an asymptotic type identity is probably the
best we can hope for in this case. Finally, to illustrate some of
the extremal cases covered by our results, we give some examples
in \S \ref{sec:geometric}.

\begin{rmk}
McShane's original identity (\ref{eqn:mcshanes id}) has also been
generalized to more general hyperbolic surfaces with cusps by
McShane himself \cite{mcshane1998im}, to hyperbolic surfaces with
cusps and/or geodesic boundary components by Mirzakhani
\cite{mirzakhani2004preprint}, to hyperbolic surfaces with cusps,
geodesic boundary and/or conical singularities, and to classical
Schottky groups by the authors in
\cite{tan-wong-zhang2004cone-surfaces},
\cite{tan-wong-zhang2004schottky}. It would be interesting to
extend the results here to representations corresponding to these
generalizations.
\end{rmk}

\vskip 5pt

\noindent {\it Acknowledgements.} We would like to thank  Makoto
Sakuma for his encouragement and for bringing our attention to
Proposition 5.2 of \cite{akiyoshi-miyachi-sakuma2004cm355}. We
would like to apologize to the authors of of
\cite{akiyoshi-miyachi-sakuma2004cm355} for the oversight in a
previous version of this paper.


\vskip 15pt
\section{{\bf Generalized Markoff Maps}}\label{sec:Markoffmaps}
\vskip 10pt

In this section we reformulate Theorems \ref{thm:A}, \ref{thm:B}
and \ref{thm:C} in terms of (generalized) Markoff maps, following
\cite{bowditch1998plms}, \cite{bowditch1997t} and
\cite{tan-wong-zhang2004gMm}, and also state our last result,
Theorem \ref{thm:D}. The discussion will be brief. The reader is
referred to the references listed above for details.

\vskip 5pt

\begin{figure}
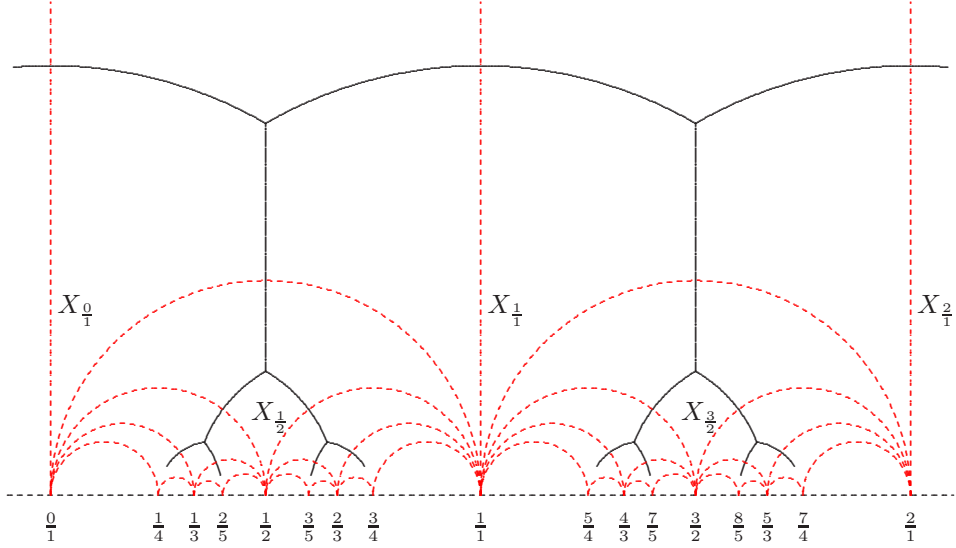


\begin{center}
\mbox{
\beginpicture
\setcoordinatesystem units <2.25in,2.25in>

\setplotarea x from -0.1 to 2.1, y from 0 to 1.2


\plot 0.5 0.86603 0.5 0.28868 /

\circulararc 38.21321 degrees from 0.5 0.28868 center at 0.66667 0

\circulararc -38.21321 degrees from 0.5 0.28868 center at 0.33333
0


 \circulararc 49.58256 degrees from 0.35714 0.12372  center at
0.375 0

\circulararc -25.03966 degrees from 0.35714 0.12372 center at 0.2
0

\circulararc 25.03966 degrees from 0.64286 0.12372 center at 0.8 0

\circulararc -49.58256 degrees from 0.64286 0.12372  center at
0.625 0


\plot 1.5 0.86603 1.5 0.28868 /

\circulararc 38.21321 degrees from 1.5 0.28868  center at 1.66667
0

\circulararc -38.21321 degrees from 1.5 0.28868 center at 1.33333
0

\circulararc 49.58256 degrees from 1.35714 0.12372  center at
1.375 0

\circulararc -25.03966 degrees from 1.35714 0.12372 center at 1.2
0

\circulararc 25.03966 degrees from 1.64286 0.12372 center at 1.8 0

\circulararc -49.58256 degrees from 1.64286 0.12372  center at
1.625 0


\circulararc 60 degrees from 1.5 0.86603  center at 1 0

\circulararc 35 degrees from 0.5 0.86603  center at 0 0

\circulararc -35 degrees from 1.5 0.86603  center at 2 0

\put {\mbox{\small $\frac{0}{1}$}} [cb] <0mm,-6mm> at 0 0

\put {\mbox{\small $\frac{1}{1}$}} [cb] <0mm,-6mm> at 1 0

\put {\mbox{\small $\frac{1}{2}$}} [cb] <0mm,-6mm> at 0.5 0

\put {\mbox{\small $\frac{1}{3}$}} [cb] <0mm,-6mm> at 0.33333 0

\put {\mbox{\small $\frac{2}{3}$}} [cb] <0mm,-6mm> at 0.66667 0

\put {\mbox{\small $\frac{1}{4}$}} [cb] <0mm,-6mm> at 0.25 0

\put {\mbox{\small $\frac{2}{5}$}} [cb] <0mm,-6mm> at 0.4 0

\put {\mbox{\small $\frac{3}{5}$}} [cb] <0mm,-6mm> at 0.6 0

\put {\mbox{\small $\frac{3}{4}$}} [cb] <0mm,-6mm> at 0.75 0


\put {\mbox{\small $\frac{2}{1}$}} [cb] <0mm,-6mm> at 2 0

\put {\mbox{\small $\frac{3}{2}$}} [cb] <0mm,-6mm> at 1.5 0

\put {\mbox{\small $\frac{4}{3}$}} [cb] <0mm,-6mm> at 1.33333 0

\put {\mbox{\small $\frac{5}{3}$}} [cb] <0mm,-6mm> at 1.66667 0

\put {\mbox{\small $\frac{5}{4}$}} [cb] <0mm,-6mm> at 1.25 0

\put {\mbox{\small $\frac{7}{5}$}} [cb] <0mm,-6mm> at 1.4 0

\put {\mbox{\small $\frac{8}{5}$}} [cb] <0mm,-6mm> at 1.6 0

\put {\mbox{\small $\frac{7}{4}$}} [cb] <0mm,-6mm> at 1.75 0

\put {\mbox{ $X_{\frac 01}$}} [cb] <0mm,-6mm> at 0.05 0.5

\put {\mbox{ $X_{\frac 11}$}} [cb] <0mm,-6mm> at 1.05 0.5

\put {\mbox{ $X_{\frac 21}$}} [cb] <0mm,-6mm> at 2.05 0.5

\put {\mbox{ $X_{\frac 12}$}} [cb] <0mm,-6mm> at 0.5 0.25

\put {\mbox{ $X_{\frac 32}$}} [cb] <0mm,-6mm> at 1.5 0.25


\setdashes<1.9pt>

\plot -0.1 0 2.1 0 /

{\color{red}

\plot 0 0 0 1.15 / \plot 1 0 1 1.15 /

\circulararc 180 degrees from 1 0  center at 0.5 0

\circulararc 180 degrees from 0.5 0  center at 0.25 0

\circulararc 180 degrees from 1 0  center at 0.75 0

\circulararc 180 degrees from 0.33333 0  center at 0.16667 0

\circulararc 180 degrees from 0.5 0  center at 0.41667 0

\circulararc 180 degrees from 0.66667 0  center at 0.58333 0

\circulararc 180 degrees from 1 0  center at 0.83333 0

\circulararc 180 degrees from 0.25 0  center at 0.125 0

\circulararc 180 degrees from 0.33333 0  center at 0.29167 0

\circulararc 180 degrees from 0.4 0  center at 0.36667 0

\circulararc 180 degrees from 0.5 0  center at 0.45 0

\circulararc 180 degrees from 0.6 0  center at 0.55 0

\circulararc 180 degrees from 0.66667 0  center at 0.63333 0

\circulararc 180 degrees from 0.75 0  center at 0.70833 0

\circulararc 180 degrees from 1 0  center at 0.875 0

\plot 2 0 2 1.15 /

\circulararc 180 degrees from 2 0  center at 1.5 0

\circulararc 180 degrees from 1.5 0  center at 1.25 0

\circulararc 180 degrees from 2 0  center at 1.75 0

\circulararc 180 degrees from 1.33333 0  center at 1.16667 0

\circulararc 180 degrees from 1.5 0  center at 1.41667 0

\circulararc 180 degrees from 1.66667 0  center at 1.58333 0

\circulararc 180 degrees from 2 0  center at 1.83333 0

\circulararc 180 degrees from 1.25 0  center at 1.125 0

\circulararc 180 degrees from 1.33333 0  center at 1.29167 0

\circulararc 180 degrees from 1.4 0  center at 1.36667 0

\circulararc 180 degrees from 1.5 0  center at 1.45 0

\circulararc 180 degrees from 1.6 0  center at 1.55 0

\circulararc 180 degrees from 1.66667 0  center at 1.63333 0

\circulararc 180 degrees from 1.75 0  center at 1.70833 0

\circulararc 180 degrees from 2 0  center at 1.875 0

}
\endpicture
}\end{center}

\caption{Farey tessellation and  the binary tree $\Sigma$.
}\label{fig:Farey}
\end{figure}

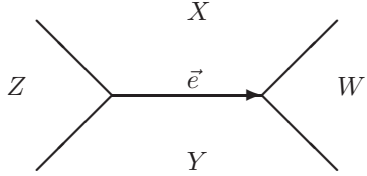
\begin{figure}
\setlength{\unitlength}{1mm}
\begin{picture}(60,30)
\thicklines \put(20,10){\vector(1,0){20}}
\put(20,10){\line(-1,1){10}} \put(20,10){\line(-1,-1){10}}
\put(40,10){\line(1,1){10}} \put(40,10){\line(1,-1){10}}
\put(30,20){$X$} \put(30,0){$Y$} \put(50,10){$W$} \put(6,10){$Z$}
\put(30,11){$\vec e$}
\end{picture}
\caption{The directed edge $\vec e = (X,Y;Z \rightarrow W)$}
\label{fig:edge oriented}
\end{figure}

Let $\Sigma$ be a binary tree properly embedded in the
(hyperbolic) plane. A {\it complementary region} of $\Sigma$ is
the closure of a connected component of the complement. Denote by
$V(\Sigma)$, $E(\Sigma)$ and $\Omega(\Sigma)$ the set of vertices,
edges, and complementary regions of $\Sigma$ respectively, or
simply, just $V$, $E$ and $\Omega$ respectively. We will fix a
concrete realization of the above concepts by thinking of $\Sigma$
as the dual to the Farey tessellation ${\mathcal{F}}$ of the
hyperbolic plane (upper half plane) by ideal triangles. Recall
that $\frac pq, \frac {p'}{q'} \in {\mathbf Q}\cup \{\infty\}$
(here $\infty:=\frac10$) with ${\rm gcd}(p,q)={\rm gcd}(p',q')=1$
are Farey neighbors if $|pq'-p'q|=1$, and that the Farey
triangulation consists of edges which are complete hyperbolic
geodesics joining all pairs $\{\frac pq, \frac {p'}{q'}\} \subset
{\mathbf Q} \cup \{\infty\}$ which are Farey neighbors, see Figure
\ref{fig:Farey}. In this way, there is a natural action of ${\rm
PSL}(2, \mathbf Z)$ on $\Sigma$, also on $\Omega$, and there is a
natural correspondence of $\Omega$ with ${\mathbf
Q}\cup\{\infty\}$, which embeds as a dense subset of ${\mathbf
R}\cup \{\infty\} \cong S^1$, and inherits the cyclic ordering
from $S^1$ ($S^1$ can be identified with the projective lamination
space of ${\mathbb T}$). Indeed, with this identification, the
action of $\PSLTwoZ$ on $\Omega \leftrightarrow {\mathbf
Q}\cup\{\infty\}$ is the usual one. We use the letters $X,Y,Z,W,
\ldots$ to denote the elements of $\Omega$, and also introduce the
notation $X_{\frac pq}$ to indicate that $X_{\frac pq} \in \Omega$
corresponds to $\frac pq \in \mathbf Q \cup \{\infty\}$, where we
use $\frac 10$ to denote $\infty$, and $\frac pq$ (with $p \in
\mathbf Z$, $q \in \mathbf N$, ${\rm gcd}(p,q)=1$) to denote
elements of $\mathbf Q$. We use the notation $e \leftrightarrow
(X,Y;Z,W)$ to indicate that $e=X \cap Y$ and $e\cap Z$ and $e \cap
W$ are the endpoints of $e$. Denote by $\vec E(\Sigma)$ (or just
$\vec E$) the set of directed edges of $\Sigma$ where the
direction is always taken to be from from the tail to the head (as
in the direction of the arrow). We also use $e$ to denote the
underlying undirected edge corresponding to the directed edge
$\vec e$, and $-\vec e$ to denote the directed edge in the
opposite direction of $\vec e$. And for a directed edge $\vec e
\in \vec E(\Sigma)$, we use $\vec e = (X,Y;Z\rightarrow W)$ to
indicate that $e \cap W$ is the head of $\vec e$, that is, $\vec
e$ is the directed edge from $Z$ to $W$, as shown in Figure
\ref{fig:edge oriented}.

\vskip 5pt

\begin{defn}\label{defn:mu-Markoff triple}
For a complex number $\mu \in \mathbf C$, a $\mu$-{\it Markoff
triple} is an ordered triple $(x,y,z) \in \mathbf C^3$ of complex
numbers satisfying the $\mu$-Markoff equation:
\begin{eqnarray}\label{eqn:x^2+y^2+z^2-xyz=mu}
x^2+y^2+z^2-xyz=\mu.
\end{eqnarray}
\end{defn}
Thus Markoff triples are just $0$-Markoff triples. It is easily
verified that if $(x,y,z)$ is a $\mu$-Markoff triple, so are
$(x,y,xy-z)$, $(x,xz-y,z)$, $(yz-x,y,z)$ and the permutation
triples of each of them.

\vskip 5pt

\begin{defn}
A map $\phi : \Omega \rightarrow \mathbf C$ is a $\mu$-{\it
Markoff map} if \begin{itemize}

\item[(i)] for every vertex $v \in V(\Sigma)$, the triple
$(\phi(X), \phi(Y), \phi(Z))$ is a $\mu$-Markoff triple, where
$X,Y,Z \in \Omega$ are the three regions meeting $v$; and

\item[(ii)] for every edge $e \in E(\Sigma)$ such that $e \leftrightarrow (X,Y;Z,W)$, we have
\begin{eqnarray}\label{eqn:xy=z+w}
xy=z+w,
\end{eqnarray}
where $x=\phi(X), y=\phi(Y)$, $z=\phi(Z)$,
$w=\phi(W)$.\end{itemize}
\end{defn}

As in \cite{tan-wong-zhang2004gMm} we shall use ${\bf \Phi}_{\mu}$
to denote the set of all $\mu$-Markoff maps, so that the set of
all Markoff maps is  denoted by ${\bf \Phi}_{0}$ in our notation.
We shall also use the following convention of Bowditch
\cite{bowditch1998plms}.

\vskip 5pt

 {\bf Convention.}\,\, We always use upper case letters
denote the regions and the corresponding lower case letters denote
the values of $\phi$ on the regions, that is, $\phi(X)=x$,
$\phi(Y)=y$, $\phi(Z)=z$, etc.

\vskip 5pt

As in the case of Markoff maps, if the edge relation
(\ref{eqn:xy=z+w}) is satisfied along all edges, then it suffices
that the vertex relation (\ref{eqn:x^2+y^2+z^2-xyz=mu}) be
satisfied at a single vertex. In fact one may establish a
bijective correspondence between $\mu$-Markoff maps and
$\mu$-Markoff triples, by fixing three regions $X, Y, Z$  which
meet at some vertex $v$; for concreteness, let us say $X=X_{\frac
01}, Y=X_{\frac 11}$ and $Z=X_{\frac 10}$. The correspondence is
then given by $\iota:\phi \mapsto (\phi(X_{\frac 01}),
\phi(X_{\frac 11}), \phi(X_{\frac 10}))$. This process may be
inverted to obtain $\iota^{-1}$ by constructing a tree of
$\mu$-Markoff triples as Bowditch did in \cite{bowditch1998plms}
for Markoff triples: Given a $\mu$-Markoff triple $(x,y,z)$, we
set $\phi(X_{\frac 01})=x, \phi(X_{\frac 11})=y, \phi(X_{\frac
10})=z$, and extend over $\Omega$ as dictated by the edge relation
(\ref{eqn:xy=z+w}). In this way one obtains an identification of
${\bf \Phi}_{\mu}$ with the variety in $\mathbf C^3$ given by the
$\mu$-Markoff equation. In particular, ${\bf \Phi}_{\mu}$ gets a
nice topology as a subset of $\mathbf C^3$.

\vskip 5pt

The natural action of  ${\rm PSL}(2, \mathbf Z)$ on $\Sigma$ and
on $\Omega$ induces an action on ${\bf \Phi}_{\mu}$, where if $H
\in \PSLTwoZ$ then $H:{\bf \Phi}_{\mu} \rightarrow {\bf
\Phi}_{\mu}$ is given by
$$H(\phi)(X)=\phi(H(X)), ~~~X \in \Omega.$$ There is also an action
of the Klein-four group, $\mathbf Z_{2}^{2}$, on ${\bf
\Phi}_{\mu}$ which reverses the signs of two of the values in each
triple,  and relates representations into $\SLTwoC$ with
representations into $\PSLTwoC$; however, this will not concern us
here.

\vskip 5pt

The key point to obtain the reformulation is that there is also a
natural correspondence between conjugacy classes of
$(\mu-2)$-representations $\rho \in {\mathcal X}_{\mu-2}$ and
$\mu$-Markoff maps, once the four regions $X,Y,Z,W$ corresponding
to a fixed edge  $e \leftrightarrow (X,Y;Z,W)$ are identified with
generators $a,b,ab$ and $ab^{-1}$ for $\Gamma$ respectively. This
works because the relations (\ref{eqn:x^2+y^2+z^2-xyz=mu}) and
(\ref{eqn:xy=z+w}) are equivalent to the trace relations. This
also gives a natural identification of $\Omega$ with $\hat \Omega$
(recall that $\hat \Omega \subset \Gamma/\sim$ is the subset
corresponding to the non-trivial, non-peripheral simple closed
curves on ${\mathbb T}$) and the action of the mapping class group
$\MCG \cong \SLTwoZ$ on ${\mathcal X}_{\mu-2}$ descends to the
$\PSLTwoZ$ action on ${\bf \Phi}_{\mu}$ described earlier.
Theorems \ref{thm:A}, \ref{thm:B} and \ref{thm:C} can then be
reformulated in terms of $\mu$-Markoff maps as follows (we use the
same definitions for the functions $h$ and ${\mathfrak h}_{\tau}$
as in the original statements).

\begin{thm}{\rm (Reformulation of Theorem \ref{thm:A})}\label{thm:AA}
Suppose $\phi \in {\bf \Phi}_{0}$. Then
\begin{eqnarray}\label{eqn:Bowditch2}
\sum_{X \in  \Omega}h\big(\phi(X))\big)=\frac{1}{2},
\end{eqnarray}
and the sum converges absolutely, if and only if

{\rm (i$'$)}\, $\phi(X) \not\in (-2,2)$ for all $X \in \Omega$,
and

{\rm (ii)} $|\phi(X)| \le 2$ for only finitely many $X \in
\Omega$.
\end{thm}

\begin{thm}{\rm (Reformulation of Theorem \ref{thm:B})}\label{thm:BB}
Suppose $\phi \in {\bf \Phi}_{\tau+2}$ {\rm(}$\tau \neq \pm
2${\rm)} and let $\nu=\cosh^{-1}(-\tau/2)$. Then
\begin{eqnarray} \label{eqn:TWZ2}
\sum_{X \in \Omega}{\mathfrak h}_{\tau} \big(\phi(X)\big)=\nu \mod
2\pi i,
\end{eqnarray}
and the sum converges absolutely, if and only if

{\rm (i$'$)} $\phi(X) \not\in (-2,2)$ for all $X \in \Omega$, and

{\rm (ii)} $|\phi(X)| \le 2$ for only finitely many $X \in
\Omega$.
\end{thm}

\begin{thm}{\rm (Reformulation of Theorem \ref{thm:C})}\label{thm:CC}
Suppose that $\phi \in {\bf \Phi}_{\tau+2}$ {\rm(}$\tau \neq
2${\rm)} is stabilized by $\langle H \rangle <\PSLTwoZ$, where $H$
is a parabolic element of $\PSLTwoZ$. Let $X_0$ be the unique
element in $ \Omega$ fixed by $H$. Then
\begin{eqnarray}\label{eqn:parabolicstab2}
\sum_{[X] \in \OmegaH-[X_0]}{\mathfrak h}\big(\phi[X]\big)=0 \mod 2\pi i, %
\end{eqnarray}
{\rm(}where ${\mathfrak h}=h$ if $\tau=-2$ and ${\mathfrak
h}={\mathfrak h}_{\tau}$ if $\tau \neq -2${\rm)} and the sum
converges absolutely, if and only if

{\rm (i)}\,\, $\phi[X] \not\in (-2,2)$ for all $[X] \in \OmegaH -
[X_0],$ and

{\rm (ii)}  $|\phi[X]| \le 2$ for only finitely many classes $[X]
\in {\OmegaH}$.
\end{thm}

Note that in Theorem \ref{thm:CC} above, we have $\phi(X_0)=2\cos
\big(\frac{q}{p}\pi\big)$ for some $\frac{q}{p} \in {\mathbf Q}$.
See Lemma \ref{lem:x0=2cos} for a proof. Now suppose that
$\phi(X_0) \in (-2,2)$ and $\phi(X_0) \neq 2\cos
\big(\frac{q}{p}\pi\big)$ for any $\frac{q}{p} \in {\mathbf Q}$.
Suppose further that for each directed edge $\vec \varepsilon$
meeting $X_0$ only at its head, $\phi$ satisfies the extended
BQ-conditions on $\Omega^{0-}(\vec \varepsilon)$ (see Definition
\ref{def:BQmaps} below). Then we have, in a special case, the
following asymptotic averaging and ergodic version of Theorem
\ref{thm:CC}.

\begin{thm}\label{thm:D}
Suppose $\phi \in {\bf \Phi}_\mu$ {\rm(}$\mu \neq 4${\rm)} is such
that, under a specific identification of $\Omega$ with ${\mathbf
Q}\cap\{\frac{1}{0}\,\}$, $\Omega_\phi(2)=\{X_{\frac10}\}$ and
$\phi(X_{\frac10}) \in (-2,2)$. Then
\begin{eqnarray}\label{eqn:phi(X0)in(-2,2)averaging}
\lim_{N \rightarrow \infty} \Re \bigg(\frac{\sum_{r \in {\mathbf Q}\cap(0,N]}{\mathfrak h}(\phi(X_r))}{N}\bigg)=0,%
\end{eqnarray}
and there exists an increasing sequence $(N_k)_{k=1}^{\infty}$ of
positive integers such that
\begin{eqnarray}\label{eqn:phi(X0)in(-2,2)asymptotic}
\lim_{k \rightarrow \infty} \textstyle{\sum_{r \in {\mathbf Q}\cap(0,N_k]}}\,{\mathfrak h}(\phi(X_r))=0 \mod 2 \pi i,%
\end{eqnarray}
where ${\mathfrak h}=h$ if $\mu=0$, and ${\mathfrak h}={\mathfrak
h}_{\mu-2}$ if $\mu \neq 0$.
\end{thm}

\vskip 3pt

We end this section by giving the formal definition and notation
for generalized Markoff maps satisfying the BQ-conditions or the
extended BQ-conditions.

\begin{defn}\label{def:BQmaps}
First, we denote by $({\bf \Phi}_{\mu})_Q$ the subset of ${\bf
\Phi}_{\mu}$ satisfying {\it the BQ-conditions}, that is, $\phi
\in({\bf \Phi}_{\mu})_Q$ if

{\rm (i)}\, $\phi(X) \not\in [-2,2]$ for all $X \in \Omega$, and

{\rm (ii)} $|\phi(X)| \le 2$ for only finitely many $X \in
\Omega$.

\vskip 3pt

Next, we denote by $({\bf \Phi}_{\mu})_{\overline Q}$ the subset
of ${\bf \Phi}_{\mu}$ satisfying {\it the extended BQ-conditions},
that is, $\phi \in({\bf \Phi}_{\mu})_{\overline Q}$ if

{\rm (i$'$)} $\phi(X) \not\in (-2,2)$ for all $X \in \Omega$, and

{\rm (ii)} $|\phi(X)| \le 2$ for only finitely many $X \in
\Omega$.

\vskip 3pt

Finally, we say that $\phi$ satisfies the BQ-conditions (resp. the
extended BQ-conditions) on a subset $\Omega'$ of $\Omega$ if (i)
(resp. (i$'$)) and (ii) above hold for all $X \in \Omega'$.
\end{defn}


 \vskip 10pt
\section{\bf Proofs of Theorems}\label{sec:proofs}
 \vskip 10pt

We first give a quick sketch of how the corresponding results were
obtained for (generalized) Markoff maps $\phi$ satisfying the
BQ-conditions. We need to recall some definitions from
\cite{bowditch1998plms} and \cite{tan-wong-zhang2004gMm}.

\vskip 5pt

Recall that $E=E(\Sigma)$ and $\vec E=\vec E(\Sigma)$ are the set
of edges and directed edges of $\Sigma$ respectively, where the
direction for a directed edge is taken to be from the tail to the
head. For $e \leftrightarrow (X,Y;Z,W)\in E$, define
$\Omega^0(e)=\{X,Y\}$. For $\vec e \leftrightarrow (X,Y;Z
\rightarrow W)\in \vec E$, $\Sigma \setminus \mbox{int}(e)$
consists of two components, denoted by $\Sigma^+(\vec e)$ and
$\Sigma^-(\vec e)$, where $\Sigma^+(\vec e)$ is the component
containing the head of $\vec e$ and $\Sigma^-(\vec e)$ is the
component containing the tail of $\vec e$ (so that $Z$ meets
$\Sigma^-(\vec e)$ and $W$ meets $\Sigma^+(\vec e)$). Define
$\Omega^{+}(\vec e)$ (resp. $\Omega^{-}(\vec e)$) to be the set of
regions in $\Omega$ whose boundaries lie entirely in
$\Sigma^{+}(\vec e)$ (resp. $\Sigma^{-}(\vec e)$). Hence $\Omega =
\Omega^{+}(\vec e) \sqcup \Omega^0(e) \sqcup \Omega^{-}(\vec e)$.
Define $\Omega^{0+}(\vec e)=\Omega^0(e) \cup \Omega^{+}(\vec e)$
(resp. $\Omega^{0-}(\vec e)=\Omega^0(e) \cup \Omega^{-}(\vec e)$).

\vskip 5pt

In terms of a concrete identification of $\Omega$ with ${\mathbf
Q}\cup\{\infty\}$, we have the following interpretation of the
above definitions. If $\vec e =
(X_{r_1},X_{r_2};X_{r_3}\rightarrow X_{r_4})\in \vec E$, then
$\Omega^0(e)\leftrightarrow \{r_1, r_2\}$ (note that $r_1$ and $
r_2$ are Farey neighbors). Also, $r_1, r_2$ divides ${\mathbf
R}\cup \{\infty\} \cong S^1$ into two open intervals $I^+(\vec e)$
and $I^-(\vec e)$ where $r_3 \in I^-(\vec e)$ and $r_4 \in
I^+(\vec e)$. (Bowditch \cite{bowditch1998plms} calls the closure
of $I^-(\vec e)$ the impression of the branch $\Omega^-(\vec e)$.)
Then $\Omega^-(\vec e)\leftrightarrow ({\mathbf Q} \cup
\{\infty\})\cap I^-(\vec e) $. For example, in Figure
\ref{fig:Farey}, if $\vec e = (X_{\frac 12},X_{\frac
11};X_{\frac23}\rightarrow X_{\frac 01})$, then $\Omega^0(\vec
e)=\{X_{\frac 12}, X_{\frac 11}\}$, and $\Omega^-(\vec
e)=\{X_{\frac pq} :\frac 12<\frac pq <\frac 11\}$. In other words,
to each directed edge $\vec e$, there corresponds an open interval
$I^-(\vec e)$ of $S^1$ with rational endpoints which are Farey
neighbors. This can be visualized as the tail of a ``comet''
(where one thinks of the directed edge as a comet); $\Omega^-(\vec
e)$ is then the ``rational dust'' in this tail.

\begin{figure}
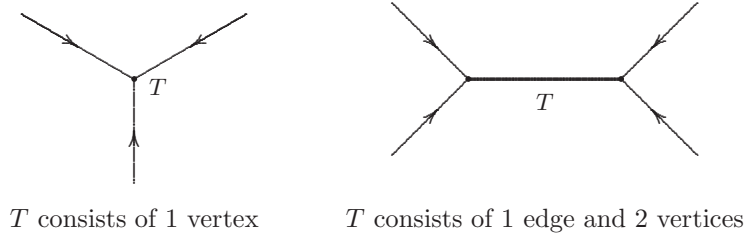


\begin{center}

\mbox{\beginpicture \setcoordinatesystem units <0.68in,0.68in>

\setplotarea x from -1.2 to 1.2, y from -1.2 to 1.2

\plot 0 -0.8 0 0 -0.866 0.5  /

\plot 0 0 0.866 0.5 /

\arrow <6pt> [.16,.6] from 0 -0.5 to 0 -0.4

\arrow <6pt> [.16,.6] from -0.866 0.5  to -0.433 0.25

\arrow <6pt> [.16,.6] from 0.866 0.5  to 0.433 0.25

\put {\mbox{\Huge $\cdot$}} [cc] <0mm,-0.2mm> at 0 0

\put {\mbox{$T$ consists of 1 vertex}} [ct] <0mm,0mm> at 0 -1.0

\put {\mbox{\small $T$}} [lc] <2mm,-1mm> at 0 0

\endpicture} \hspace{0.2in} \raisebox{-11ex}

\mbox{\beginpicture \setcoordinatesystem units <0.4in,0.4in>

\setplotarea x from -1.2 to 3.2, y from -1.2 to 1.2

\plot -1 1 0 0 2 0 3 1 /

\plot -1 -1 0 0 / \plot 0 0.01 2 0.01 /

\plot 0 -0.01 2 -0.01 / \plot 2 0 3 -1 /

\arrow <6pt> [.16,.6] from  -0.5 0.5  to -0.4 0.4

\arrow <6pt> [.16,.6] from -0.5 -0.5 to -0.4 -0.4

\arrow <6pt> [.16,.6] from 2.5 0.5  to 2.4 0.4

\arrow <6pt> [.16,.6] from 2.5 -0.5  to 2.4 -0.4

\put {\mbox{\Huge $\cdot$}} [cc] <0mm,-0.2mm> at 0 0

\put {\mbox{\Huge $\cdot$}} [cc] <0mm,-0.2mm> at 2 0

\put {\mbox{$T$ consists of 1 edge and 2 vertices}} [ct] <0mm,0mm>
at 1 -1.7

\put {\mbox{\small $T$}} [cc] <0mm,-3mm> at 1 0

\endpicture}
\end{center}
\caption{Two simple circular sets $C(T)$}\label{fig:circular}
\end{figure}

\vskip 5pt

For a given subtree $T \subseteq \Sigma$, the set $C(T)\subseteq
\vec E$ is the subset of directed edges $\vec e$ such that $\vec e
\in C(T)$ if and only if $e \cap T$ consists of only one point,
the point being the head of $\vec e$. If $T$ is finite, then so is
$C(T)$. For example, if $T$ consists of a single vertex $v$, then
$C(T)$ consists of the three directed edges with heads at $v$, and
if $T$ consists of a single edge and its two vertices, then $C(T)$
consists of four directed edges, as shown in Figure
\ref{fig:circular}. A finite subset $C$ of $\vec E$ is said to be
a {\it circular set} if it has the form $C(T)$ for some finite
subtree $T$ of $\Sigma$. It is not difficult to see, using the
concrete realization, that $C$ is either circular or of the form
$\{\vec e, -\vec e\}$ for some $\vec e \in \vec E$\, if and only
if
\begin{eqnarray*}
\displaystyle \bigcup_{\vec e \in C} \overline{I^-(\vec
e)}&=&{\mathbf R}\cup \{\infty\} \cong S^1, \quad \hbox{
and}\\
I^-(\vec e_1)\cap I^-(\vec e_2)&=&\emptyset, ~~~
\,\,\,\hbox{for}~~~\,\, \vec e_1, \vec e_2 \in C, \,\, \vec e_1
\neq \vec e_2.
\end{eqnarray*}
Using the analogy of comets as directed edges, this corresponds to
a finite set of comets whose tails fill the entire horizon at
infinity, and such that no comet lies in the ``shadow'' of
another.

\vskip 5pt

Now let us fix one $\phi \in (\Phi_{\mu})_Q$. (We shall deal with
Theorems \ref{thm:BowditchT} and \ref{thm:CC} later.) So $\phi$
satisfies the BQ-conditions (see Definition \ref{def:BQmaps}).
Using $\phi$ we can define a map
$$\alpha=\alpha_\phi: E\rightarrow \vec E$$ such that, for $e\leftrightarrow(X,Y;Z,W)$, $$\alpha(e)
= (X,Y;Z\rightarrow W)$$ if $|\phi(Z)| > |\phi(W)|$ and one may
choose an arbitrary direction if $|\phi(Z)|=|\phi(W)|$. Thus the
direction of the edge is always from the region with greater
absolute value of $\phi$ to the region with smaller absolute value
of $\phi$.

There is a way to measure the growth rate of the function
$\log^+|\phi|$ on $\Omega$ (where $\log^+ x:=\mbox{max}\{0, \log
x\}$), which will be useful for proving absolute convergence of
certain series, by comparing it to the Fibonacci function
$F:\Omega \rightarrow {\mathbf N}$, which we can take to be
$$F(X_{\frac pq})=|p|+|q|, \,\,\,\mbox{where}\,\,\,\, {\rm gcd}(p,q)=1.$$
(Strictly speaking, this is the Fibonacci function associated to
the edge $e\leftrightarrow (X_{\frac10},$
$X_{\frac01};X_{\frac{-1}1}, X_{\frac1{1}})$; but it was shown in
\cite{bowditch1998plms} that in fact the concept of having
Fibonacci bounds is independent of the edge used.)  Then
$\log^+|\phi|$ has a lower Fibonacci bound on $\Omega$ if there
exists a constant $\kappa > 0$ such that $\log^+|\phi(X)| > \kappa
\, F(X)$ for all but finitely many $X \in \Omega$. Clearly, we can
talk about lower Fibonacci bounds for $\log^+|\phi|$ restricted to
subsets $\Omega'\subseteq \Omega$ as well. The key point is that
if $\log^+|\phi|$ has a lower Fibonacci bound on $\Omega^-(\vec
e)$ for some $\vec e \in \vec E$, then the sum $\sum_{X \in
\Omega^-(\vec e)} {\mathfrak h}_{\tau}(\phi(X))$ converges
absolutely (to a value which can be determined).

\vskip 5pt

We now outline (with some slight modifications) in several steps
the basic strategy used in \cite{bowditch1998plms} and
\cite{tan-wong-zhang2004gMm} to prove the (generalized) McShane's
identity for a fixed $\phi \in (\Phi_{\mu})_Q$. For any $k>0$, we
define $$\Omega(k):= \Omega_{\phi}(k)= \{X \in \Omega \mid
|\phi(X)| \le k\}.$$

By BQ(ii), $\Omega(2)$ is finite. Let us first suppose that
$\Omega(2)$ is non-empty.

\vskip 5pt


\begin{figure}
\begin{center}\mbox{
\hskip -0.05in
\beginpicture

\setcoordinatesystem units <0.55in,0.55in> \setplotarea x from
-0.5 to 7.5, y from 0.3 to 2.8

\plot 0 1 7 1 /

\arrow <6pt> [.16,.6] from 0.5 1  to 0.6 1

\arrow <6pt> [.16,.6] from 1.5 1  to 1.6 1

\arrow <6pt> [.16,.6] from 2.5 1  to 2.6 1

\arrow <6pt> [.16,.6] from 3.5 1  to 3.6 1

\arrow <6pt> [.16,.6] from 4.5 1  to 4.6 1

\arrow <6pt> [.16,.6] from 5.5 1  to 5.6 1

\arrow <6pt> [.16,.6] from 6.5 1  to 6.6 1

\plot 1 2 1 1 /  \plot 2 2 2 1 /

\plot 3 2 3 1 /  \plot 4 2 4 1 /

\plot 5 2 5 1 /  \plot 6 2 6 1 /

\plot 0.7 2.3 1 2 1.3 2.3 /

\plot 1.7 2.3 2 2 2.3 2.3 /

\plot 2.7 2.3 3 2 3.3 2.3 /

\plot 3.7 2.3 4 2 4.3 2.3 /

\plot 4.7 2.3 5 2 5.3 2.3 /

\plot 5.7 2.3 6 2 6.3 2.3 /

\put {\mbox{\Huge $\cdot$}} [cc] <0mm,-0.2mm> at 1 1

\put {\mbox{\Huge $\cdot$}} [cc] <0mm,-0.2mm> at 2 1

\put {\mbox{\Huge $\cdot$}} [cc] <0mm,-0.2mm> at 3 1

\put {\mbox{\Huge $\cdot$}} [cc] <0mm,-0.2mm> at 4 1

\put {\mbox{\Huge $\cdot$}} [cc] <0mm,-0.2mm> at 5 1

\put {\mbox{\Huge $\cdot$}} [cc] <0mm,-0.2mm> at 6 1

\put {\mbox{\small ${\vec e}_{-N}$}} [ct] <0mm,-2mm> at 1.5 1

\put {\mbox{\small ${\vec e}_{N}$}} [ct] <0mm,-2mm> at 5.5 1

\put {\mbox{\small $X_0$}} [cc] <0mm,1mm> at 3.5 0.30

\put {\mbox{$\cdots$}} [cc] <-1mm,1mm> at 0.6 1.25

\put {\mbox{$\cdots$}} [cc] <-1mm,1mm> at 6.6 1.25

\put {\mbox{\small $Y_{-N}$}} [cb] <0mm,1mm> at 1.5 1.15

\put {\mbox{$\cdots$}} [cc] <-1mm,1mm> at 2.6 1.25

\put {\mbox{\small $Y_0$}} [cb] <0mm,1mm> at 3.5 1.15

\put {\mbox{$\cdots$}} [cc] <-1mm,1mm> at 4.6 1.25

\put {\mbox{\small $Y_N$}} [cb] <0mm,1mm> at 5.5 1.15

\endpicture}\end{center}

\caption{}\label{fig:aroundX}

\end{figure}


\nnn \underline{Step 1.} \\ \nnn Let $X_0 \in \Omega(2)$ and let
$\{\vec e_n \mid n\in {\mathbf Z}\}$ be the bi-infinite sequence
of directed edges in $\Sigma$ bounding $X_0$, say in clockwise
order around $X_0$, where $\vec e_n$ is directed from $\vec
e_{n-1}$ to $\vec e_{n+1}$, and let $Y_n$ be the region sharing
the edge $e_n$ with $X_0$ (see Figure \ref{fig:aroundX} where we
have represented the boundary of $X_0$ by a straight line so all
the edges $e_n$ lie on a straight line). Then by BQ(i), $\phi(X_0)
\not\in [-2,2]$, so $|\phi(Y_n)|$ grows exponentially in $|n|$ for
$|n|$ sufficiently large (Lemma 3.3 \cite{bowditch1998plms}, see
also \cite{tan-wong-zhang2004gMm}). We can find a finite subarc
$$J(X_0)=\bigcup_{-N < i < N} e_i, \qquad {\rm where} \,\, N \,\,
\textrm{is a positive integer}$$ such that

(i) if\, $Y_j \in \Omega(2)$\, then\, $e_j \subseteq J(X_0)$;

(ii) if\, $e_k \not\subseteq J(X_0)$\, then \,$\alpha(e_k)$\, is
directed towards\, $J(X_0)$,\, that is, \,$\alpha(e_k)=\vec e_k$
for $k \le -N$ and $\alpha(e_k)=-\vec e_k$ for $k \ge N$; and

(iii) $\log^{+}|\phi|$ has a lower Fibonacci bound on both
$\Omega^{-}(-\vec e_{N})$ and $\Omega^{-}(\vec e_{-N})$.

\vskip 5pt

\nnn \underline{Step 2.}\\ \nnn Now we take the union
$$T=\bigcup_{X_0 \in \Omega(2)} J(X_0).$$ It can be shown that $T$
is always connected (so $T$ is a finite subtree of $\Sigma$) and
further that $\alpha(e)$ is directed towards $T$ for all $e
\not\subseteq T$.

\vskip 5pt

\nnn\underline{Aside 1:} \\
In the case where $\Omega(2)=\emptyset$, $T$ consists of a single
vertex $v$ which is a sink, with $\alpha(e_i)$ directed towards
$v$ for all the three edges $e_1,e_2,e_3$ incident to $v$. In this
case $\log^{+}|\phi|$ has a lower Fibonacci bound on every $\vec e
\in \vec E$.

\vskip 5pt

\nnn \underline{Step 3.}\\ \nnn Now we want to show that for all
$\vec e \in C(T)$, $\log^{+}|\phi|$ has a lower Fibonacci bound on
$\Omega^{-}(\vec e)$.  There are two possibilities. First,
$|\phi(X)|>2$ for $X \in \Omega^0(\vec e)$ and furthermore
$\alpha(e)=\vec e$, then $\log^{+}|\phi|$ has a lower Fibonacci
bound on $\Omega^{-}(\vec e)$ (see Corollary 3.6
\cite{bowditch1998plms}, see also \cite{tan-wong-zhang2004gMm}).
The second possibility is that $|\phi(X)|\le 2$ for exactly one of
the $X \in \Omega^0(\vec e)$. This can only occur in the situation
of part (iii) of Step 1, with $\vec e=-\vec e_{N}$ or $\vec
e_{-N}$, so again $\log^{+}|\phi|$ has a lower Fibonacci bound on
$\Omega^-(\vec e)$.

Since $\Omega=\bigcup_{\vec e \in C(T)} \Omega^{0-}(\vec e)$, this
gives the absolute convergence of the sums in
(\ref{eqn:Bowditch2}) and (\ref{eqn:TWZ2}), as the fact that
$\log^{+}|\phi|$ has a lower Fibonacci bound implies the absolute
convergence of the sums.

\vskip 5pt

\nnn\underline{Aside 2:}\\ A slight modification is necessary to
deal with Theorems \ref{thm:BowditchT} and \ref{thm:CC}. We
construct the tree $T$ as before in Step 2, except that now it is
not finite. However, the tree is invariant under the action of
$\langle H\rangle$, and by taking the quotient of $T$ by the
action of $\langle H \rangle$, we do indeed get a finite graph.
Following \cite{bowditch1997t}, by considering the sum over
suitable branches, we again get the absolute convergence of the
sum in (\ref{eqn:hyperbolicstab}) and (\ref{eqn:parabolicstab2}),
when $\phi/\langle H\rangle$ satisfies the BQ-conditions.

\vskip 5pt

\nnn \underline{Step 4.}\\ \nnn To get the actual values of the
sums requires the introduction of a function
$$\Psi_{\mu}: \{{\mu}\textrm{-Markoff triples}\} \rightarrow
{\mathbf C}$$ (where $\mu=\tau+2\neq 0,4$) defined by
\begin{eqnarray}\label{Psi(x,y,z)=log}
\Psi_{\mu}(x,y,z):= \log
\frac{1+(e^{\nu}-1)(z/xy)}{\sqrt{1-\mu/x^2}\sqrt{1-\mu/y^2}},
\end{eqnarray}
where $\nu=\cosh^{-1}(1-\mu/2)$, for a $\mu$-Markoff triple
$(x,y,z)$ with $x,y \neq 0,\pm \sqrt{\mu}$.


\vskip 5pt

\begin{rmk}
In the case where $\mu=0$ if we set $\mu=0$ (and hence $\nu=0$) in
(\ref{Psi(x,y,z)=log}) we would have $\Psi_{0}(x,y,z) \equiv 0$;
in this case we define a function $\Psi'_0$ by
$$\Psi'_0(x,y,z) = z/xy,$$
which is used extensively by Bowditch \cite{bowditch1998plms} and
is in fact the partial derivative of $2\,\Psi_{\mu}(x,y,z)$ with
respect to $\nu$ evaluated at $\nu=0$ (or equivalently, $\mu=0$).
\end{rmk}

\vskip 5pt

Although the function $\Psi_{\mu}$ is rather complicated, what
concerns us most will be the properties of the function and the
way it relates to the sums in (\ref{eqn:hyperbolicstab}),
(\ref{eqn:Bowditch2}), (\ref{eqn:TWZ2}), and
(\ref{eqn:parabolicstab2}).

\vskip 5pt

\begin{lem}\label{lem:Psi} {\rm(Properties of the function $\Psi_\mu$)}
Suppose $\mu \neq 0,4$. \begin{itemize}

\item[{\rm(i)}] If $(x,y,z)$ is a $\mu$-Markoff triple such that $x,y,z
\neq 0, \pm \sqrt{\mu}$ then
\begin{eqnarray}\label{sumPsi i}
\Psi_{\mu}(x,y,z)+\Psi_{\mu}(y,z,x)+\Psi_{\mu}(z,x,y)=\nu \mod 2
\pi i.
\end{eqnarray}

\item[{\rm(ii)}] If $(x,y,z)$ and $(x,y,w)$ are $\mu$-Markoff triples
such that $x,y \neq 0, \pm \sqrt{\mu}$ and $z+w=xy$ then
\begin{eqnarray}\label{sumPsi ii}
\Psi_{\mu}(x,y,z)+\Psi_{\mu}(x,y,w)=\nu \mod 2 \pi i.
\end{eqnarray}
\end{itemize}
\end{lem}

\begin{defn}\label{defn:psi(e)}
For a fixed $\phi \in \Phi_{\mu}$ ($\mu \neq 0,4$), we define the {\it
edge-weight function} $\psi$ associated to $\phi$ by
\begin{eqnarray}\label{eqn:psi(e)}
\psi(\vec e):=\psi_{\phi}(\vec e)=\Psi_{\mu}(x,y,z),
\end{eqnarray}
where $\vec e = (X,Y;W\rightarrow Z)$ (that is, $Z$ is at the head
of $\vec e$) and, by our convention, $x=\phi(X), y=\phi(Y),
z=\phi(Z)$.
\end{defn}

\begin{rmk} In the case where $\mu=0$ we define
\begin{eqnarray}\label{eqn:psi(e) mu=0}
\psi(\vec e)={z}/{xy}
\end{eqnarray}
as explained before.
\end{rmk}

Then we have the following properties, proven in
\cite{tan-wong-zhang2004gMm} with the help of Lemma \ref{lem:Psi}
above.

\begin{lem}\label{lem:psi} {\rm(Properties of the edge-weight function
$\psi$)}\,\, Suppose $\phi \in \Phi_{\mu}$ {\rm(}where $\mu \neq
0,4${\rm)} and $\phi(X) \neq 0, \pm \sqrt{\mu}$ for all $X \in
\Omega$. Then

{\rm (i)} for $\vec e \in \vec E$,
\begin{eqnarray}\label{eqn:sum psi i}
\psi(\vec e)+\psi(-\vec e)=\nu \mod 2 \pi i;
\end{eqnarray}

{\rm (ii)} for a circular set $C \subset \vec E$,
\begin{eqnarray}\label{eqn:sum psi ii}
\sum_{\vec e \in C} \psi(\vec e)=\nu \mod 2 \pi i.
\end{eqnarray}
\end{lem}

\begin{rmk} In the case where $\mu = 0$, Lemmas \ref{lem:Psi} and \ref{lem:psi} hold with
$\nu$ replaced by $1$ and with ${\rm mod}\, 2 \pi i$ removed, as
used by Bowditch in \cite{bowditch1998plms}.
\end{rmk}

\begin{defn}\label{defn:hat frak h} (The function $\hat{\mathfrak h}_{\tau}$) \,\,
For $\tau \in {\mathbf C}\backslash\{\pm 2\}$, we set $\mu=\tau+2$
and define another function $\hat{\mathfrak h}_{\tau}:{\mathbf C}
\backslash \{\pm\sqrt{\mu}\} \rightarrow {\mathbf C}$, which is a
specific half, modulo $2\pi i$, of the function ${\mathfrak
h}_{\tau}$ defined before, by
\begin{eqnarray}
\hat{\mathfrak
h}_{\tau}(x)=\log\frac{1+(e^{\nu}-1)h(x)}{\sqrt{1-\mu/x^2}},
\end{eqnarray}
where $\nu=\cosh^{-1}(-\tau/2)$ and $h$ is the function defined by
(\ref{eqn:h(x)=}).
\end{defn}

\vskip 5pt

It can be checked by direct calculation that
\begin{eqnarray}\label{eqn:frakh=2hatfrakh}
{\mathfrak h}_{\tau}=2\,\hat{\mathfrak h}_{\tau} \mod 2\pi i.
\end{eqnarray}
We also have
\begin{eqnarray}\label{eqn:hatfrakh=limPsi}
\hat{\mathfrak
h}_{\tau}(x)=\lim_{y\rightarrow\infty}{\Psi}_{\mu}(x,y,z),
\end{eqnarray}
where $\mu=\tau+2$ and $z$ is given by
\begin{eqnarray}
z=\frac{xy}{2}\bigg(1-\sqrt{1-4\bigg(\frac{1}{x^2}+\frac{1}{y^2}-\frac{\mu}{x^2y^2}\bigg)}\,\,\bigg).
\end{eqnarray}
Note that the properties (\ref{eqn:frakh=2hatfrakh}) and
(\ref{eqn:hatfrakh=limPsi}) are crucial for the following theorem
to hold, as can be seen in \cite{tan-wong-zhang2004gMm}.

\vskip 5pt

\begin{thm}\label{thm:partialsum} Suppose $\phi \in {\bf \Phi}_{\mu}$ {\rm(}$\mu=\tau+2 \neq 0,4${\rm)}
and $\vec e \in \vec E$ is such that $\log^{+}|\phi|$ has a lower
Fibonacci bound on $\Omega^{-}(\vec e)$. Then
\begin{eqnarray}\label{eqn:partialsum}
\sum_{X \in \Omega^0(e)}\hat{\mathfrak h}_{\tau}(\phi(X))+\sum_{X
\in \Omega^-(\vec e)}{\mathfrak h}_{\tau}(\phi(X))=\psi(\vec e)
\mod 2 \pi i,
\end{eqnarray} where the infinite sum converges absolutely.
\square
\end{thm}

Note that the second sum on the left hand side of
(\ref{eqn:partialsum}) is an infinite sum while the first sum has
only two summands.

\vskip 5pt

Putting everything together, we can now prove Theorems
\ref{thm:AA}, \ref{thm:BB} and \ref{thm:CC} in the case where the
BQ-conditions are satisfied. We just split the infinite sum in
each case into a finite number of infinite sums over the tails of
a circular set of directed edges as in Step 3, evaluate each using
Theorem \ref{thm:partialsum}, and then apply (ii) of Lemma
\ref{lem:psi} to get the sum. A slight variation is needed for
Theorem \ref{thm:CC}: we need to choose a suitable tree $T'$ which
is a ``fundamental domain'' for the invariant tree $T$ and we use
the fact that there are two edges $\vec {e'}, \vec {e''}$ in the
circular set of $C(T')$ with $\psi(\vec {e''})=\psi(-\vec {e'})$
because $H$ maps $\vec {e'}$ to $-\vec {e''}$ (see
\cite{bowditch1997t} for details).

\vskip 5pt

Now a careful examination of the above sketch will show that the
proof in the general case breaks down only at Step 1(iii);
specifically, only in the last statement there. This is because if
we have $\phi(X_0)=\pm 2$ for some $X_0\in \Omega$, and suppose
$Y_n, n \in {\mathbf Z}$ are the neighboring regions of $X_0$,
then $|\phi(Y_n)|$ only grows linearly in $|n|$ for large $|n|$,
and not exponentially as in the case when $\phi(X_0) \not\in
[-2,2]$. So for the directed edges $-\vec e_{N}$ and $\vec e_{-N}$
directed towards $J(X_0)$ in Step 1(iii), $\log^{+}|\phi|$ does
not have lower Fibonacci bounds on $\Omega^{-}(-\vec e_{N})$ and
$\Omega^{-}(\vec e_{-N}$). Nonetheless, we show in the following
proposition that the sums in (\ref{eqn:Bowditch2}),
(\ref{eqn:TWZ2}) and (\ref{eqn:parabolicstab2}) over
$\Omega^{0-}(-\vec e_{N})$ and $\Omega^{0-}(\vec e_{-N})$ converge
absolutely, and with appropriate adjustments, to the values
$\psi(-\vec e_{N})$ and $\psi(\vec e_{-N})$ as in the case when
$\phi(X_0) \not\in [-2,2]$.

\vskip 5pt

\begin{prop}\label{lem:convergence around pm 2}
Suppose $\phi \in \Phi_{\mu}$ {\rm(}$\mu=\tau+2 \neq 0,4${\rm)}
and $\vec e \in \vec E$ with $\Omega^{0}(e)=\{X_0,Y_0\}$ such that
$\phi(X_0)=\pm 2$ and $\Omega(2) \cap \Omega^{0-}(\vec
e)=\{X_0\}$. Then
\begin{eqnarray}\label{eqn:partial sum 2}
\sum_{X \in \Omega^{0}(e)} \hat {\mathfrak
h}_{\tau}(\phi(X))+\sum_{X \in \Omega^{-}(\vec e)}{\mathfrak
h}_{\tau}(\phi(X))=\psi(\vec e) \mod 2 \pi i
\end{eqnarray}
where the infinite sum converges absolutely.
\end{prop}

\nnn {\it Proof.}\,\, For simplicity of notation, let us write
$\hat {\mathfrak h}=\hat {\mathfrak h}_{\tau}$, ${\mathfrak
h}={\mathfrak h}_{\tau}$ in this proof.

\begin{figure}
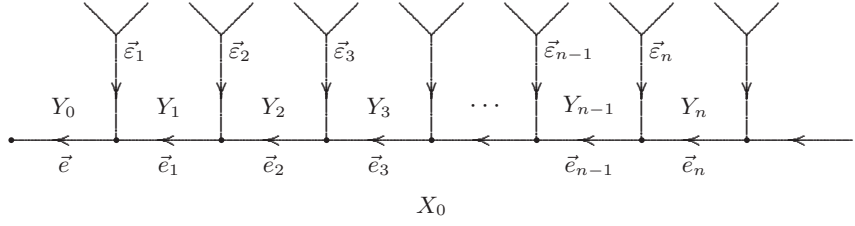

\begin{center}\mbox{
\hskip -0.05in
\beginpicture

\setcoordinatesystem units <0.55in,0.55in> \setplotarea x from
-0.5 to 8.5, y from -0.3 to 2.8

\plot 0 1 8 1 /

\arrow <6pt> [.16,.6] from 0.5 1  to 0.4 1

\arrow <6pt> [.16,.6] from 1.5 1  to 1.4 1

\arrow <6pt> [.16,.6] from 2.5 1  to 2.4 1

\arrow <6pt> [.16,.6] from 3.5 1  to 3.4 1

\arrow <6pt> [.16,.6] from 4.5 1  to 4.4 1

\arrow <6pt> [.16,.6] from 5.5 1  to 5.4 1

\arrow <6pt> [.16,.6] from 6.5 1  to 6.4 1

\arrow <6pt> [.16,.6] from 7.5 1  to 7.4 1

\plot 1 2 1 1 /  \plot 2 2 2 1 /

\plot 3 2 3 1 /  \plot 4 2 4 1 /

\plot 5 2 5 1 /  \plot 6 2 6 1 /

\plot 7 2 7 1 /

\arrow <6pt> [.16,.6] from 1 1.5 to 1 1.4

\arrow <6pt> [.16,.6] from 2 1.5 to 2 1.4

\arrow <6pt> [.16,.6] from 3 1.5 to 3 1.4

\arrow <6pt> [.16,.6] from 4 1.5 to 4 1.4

\arrow <6pt> [.16,.6] from 5 1.5 to 5 1.4

\arrow <6pt> [.16,.6] from 7 1.5 to 7 1.4

\arrow <6pt> [.16,.6] from 6 1.5 to 6 1.4

\plot 0.7 2.3 1 2 1.3 2.3 /

\plot 1.7 2.3 2 2 2.3 2.3 /

\plot 2.7 2.3 3 2 3.3 2.3 /

\plot 3.7 2.3 4 2 4.3 2.3 /

\plot 4.7 2.3 5 2 5.3 2.3 /

\plot 6.7 2.3 7 2 7.3 2.3 /

\plot 5.7 2.3 6 2 6.3 2.3 /

\put {\mbox{\Huge $\cdot$}} [cc] <0mm,-0.2mm> at 0 1

\put {\mbox{\Huge $\cdot$}} [cc] <0mm,-0.2mm> at 1 1

\put {\mbox{\Huge $\cdot$}} [cc] <0mm,-0.2mm> at 2 1

\put {\mbox{\Huge $\cdot$}} [cc] <0mm,-0.2mm> at 3 1

\put {\mbox{\Huge $\cdot$}} [cc] <0mm,-0.2mm> at 4 1

\put {\mbox{\Huge $\cdot$}} [cc] <0mm,-0.2mm> at 5 1

\put {\mbox{\Huge $\cdot$}} [cc] <0mm,-0.2mm> at 6 1

\put {\mbox{\Huge $\cdot$}} [cc] <0mm,-0.2mm> at 7 1

\put {\mbox{\small $\vec e$}} [ct] <0mm,-2mm> at 0.5 1

\put {\mbox{\small ${\vec e}_1$}} [ct] <0mm,-2mm> at 1.5 1

\put {\mbox{\small ${\vec e}_2$}} [ct] <0mm,-2mm> at 2.5 1

\put {\mbox{\small ${\vec e}_3$}} [ct] <0mm,-2mm> at 3.5 1

\put {\mbox{\small ${\vec e}_{n-1}$}} [ct] <0mm,-2mm> at 5.5 1

\put {\mbox{\small ${\vec e}_{n}$}} [ct] <0mm,-2mm> at 6.5 1

\put {\mbox{\small $X_0$}} [cc] <0mm,1mm> at 4 0.30

\put {\mbox{$\cdots$}} [cc] <-1mm,1mm> at 4.6 1.25

\put {\mbox{\small $Y_0$}} [cb] <0mm,1mm> at 0.5 1.15

\put {\mbox{\small $Y_1$}} [cb] <0mm,1mm> at 1.5 1.15

\put {\mbox{\small $Y_2$}} [cb] <0mm,1mm> at 2.5 1.15

\put {\mbox{\small $Y_3$}} [cb] <0mm,1mm> at 3.5 1.15

\put {\mbox{\small $Y_{n-1}$}} [cb] <0mm,1mm> at 5.5 1.15

\put {\mbox{\small $Y_n$}} [cb] <0mm,1mm> at 6.5 1.15

\put {\mbox{\small ${\vec \varepsilon}_1$}} [lc] <1mm,0mm> at 1
1.85

\put {\mbox{\small ${\vec \varepsilon}_2$}} [lc] <1mm,0mm> at 2
1.85

\put {\mbox{\small ${\vec \varepsilon}_3$}} [lc] <1mm,0mm> at 3
1.85

\put {\mbox{\small ${\vec \varepsilon}_{n-1}$}} [lc] <1mm,0mm> at
5 1.85

\put {\mbox{\small ${\vec \varepsilon}_{n}$}} [lc] <1mm,0mm> at 6
1.85

\endpicture}\end{center}
\vskip -0.25in \caption{$\phi(X_0)=2$}\label{fig:phi(X)=2}
\end{figure}

We may assume $\phi(X_0)=2$ since the proof for the case where
$\phi(X_0)=-2$ is essentially the same. The proof given here is
very similar to that given in
\cite{akiyoshi-miyachi-sakuma2004cm355} for the $\Phi_0$ case,
although somewhat more complicated because of the nature of
${\mathfrak h}$. Let $Y_n, n \ge 1$ be the the sequence of
neighboring regions of $X_0$ which lie in $\Omega^{-}(\vec e)$ and
let $e_n, n \ge 0$ be the edge $Y_n \cap X_0$ (so $e_0=e$). Let
$\vec e_n, n \ge 1$ be the directed edge on $e_n$ which is
directed towards $e$. Let $v_n, n \ge 1$ be the vertex between
$\vec e_{n-1}$ and $\vec e_n$. Let $\varepsilon_n, n \ge 1$ be the
edge $Y_{n-1} \cap Y_n$ and $\vec \varepsilon_n \in \vec E$ be the
directed edge on $\varepsilon_n$ with its head at $v_n$ (see
Figure \ref{fig:phi(X)=2}).

Let $y_n=\phi(Y_n), n \ge 0$. Then an easy computation using
(\ref{eqn:x^2+y^2+z^2-xyz=mu}) and (\ref{eqn:xy=z+w}) gives
$$y_n=y_0+n(\mu-4)^{1/2},$$ for some choice of the square
root. Since $\mu \neq 4$, $y_n$ grows linearly in $n$.

Since $\Omega(2) \cap \Omega^{0-}(\vec \varepsilon_n)=\emptyset$,
we know that $\log^{+}|\phi|$ has a lower Fibonacci bound on
$\Omega^{0-}(\vec \varepsilon_n)$ and hence
$\sum_{X\in\Omega^{0-}(\vec \varepsilon_n)}{\mathfrak h}(\phi(X))$
converges absolutely for all $n \ge 1$.

We need to show that $\sum_{X \in \Omega^{0-}(\vec
e)}\big|{\mathfrak h}(\phi(X))\big|$ converges. This follows from
the following lemma since
\begin{eqnarray*}
\Omega^{0-}(\vec e)=\{X_0\} \cup
\bigcup_{n=1}^{\infty}\Omega^{0-}(\vec \varepsilon_n).
\end{eqnarray*}

\begin{lem}\label{lem:estimate} There exists a constant $\kappa >
0$ such that for sufficiently large $n$
\begin{eqnarray}\label{eqn:estimate}
\sum_{X \in \Omega^{0-}(\vec \varepsilon_n)}\big|{\mathfrak
h}(\phi(X))\big| \le \kappa \cdot \frac{1}{n^2}.
\end{eqnarray}
\end{lem}

\begin{pf} There is a
constant $\kappa_1>0$ such that for all $n \ge 1$, $|y_n| \ge
\kappa_{1}n$ and hence $\log|y_n| \ge \log n + \log \kappa_1$.

Let $F_{\varepsilon_n}, n \ge 1$ be the Fibonacci function with
respect to edge $\varepsilon_n$ (see the beginning of this section
where the Fibonacci function was defined for the edge between
$X_{\frac 10}$ and $X_{\frac 01}$, see also
\cite{bowditch1998plms}). Hence there exist constants
$\kappa_2,\kappa_3, \kappa_4>0$ such that if $n$ is sufficiently
large, then
\begin{eqnarray}\label{eqn:logpsi estimate}
\log \big|\phi(X)\big| &\ge& \min\{\log |y_{n-1}|-\log 2, \log
|y_n|-\log 2\}\,F_{\varepsilon_n}(X)\\ &\ge& (\,\log n +
\log\kappa_2)\, F_{\varepsilon_n}(X) \\ &\ge& (\,\log n +
\log\kappa_3) + \kappa_4 \, F_{\varepsilon_n}(X)
\end{eqnarray}
for all $X \in \Omega^{0-}(\vec \varepsilon_n)$.
Thus there exists a constant $\kappa_5>0$ so that for sufficiently
large $n$,
\begin{eqnarray}
\textstyle{\sum_{X\in\Omega^{0-}(\vec
\varepsilon_n)}}\,\big|\phi(X)\big|^{-2} \le \kappa_5 \cdot
\frac{1}{n^2}.
\end{eqnarray}
This gives (\ref{eqn:estimate}), noticing that $\big|{\mathfrak
h}(\phi(X))\big|=O\big(\big|\phi(X)\big|^{-2}\big)$. Actually,
since $h(x) = O(x^{-2})$ as $x \rightarrow \infty$, we have
\begin{eqnarray*}
{\mathfrak h}(x)&=&\log
\bigg(\frac{1+(e^{\nu}-1)h(x)}{1+(e^{-\nu}-1)h(x)}\bigg)\\
&=&\log\bigg(1+\frac{2(\sinh\nu)h(x)}{1+(e^{-\nu}-1)h(x)}\bigg)\\
&\sim&2(\sinh\nu)h(x)\\ &=& O(x^{-2}).
\end{eqnarray*}
This proves Lemma \ref{lem:estimate}.
\end{pf}

Thus the infinite sum in (\ref{eqn:partial sum 2}) converges
absolutely. Its actual value can be obtained by sharpening the
proof of (\ref{eqn:partialsum}) in
\cite{tan-wong-zhang2004schottky}. Alternatively, we can also
evaluate it using (\ref{eqn:partialsum}) by a simple limit process
as follows.

Let $T=\bigcup_{i=1}^{n-1}e_{i}$. Applying Lemma \ref{lem:psi}(ii)
to the circular set $C(T) \subset \vec E$ gives
\begin{eqnarray*}
\textstyle{\sum_{\vec {e'} \in C(T)}}\,\psi(\vec {e'})=\nu \mod 2 \pi i. %
\end{eqnarray*}
Note that $C(T)=\{-\vec e, \vec \varepsilon_1, \cdots, \vec
\varepsilon_{n}, \vec e_n \}$. Thus by Lemma \ref{lem:psi}(i) and
Theorem \ref{thm:partialsum} we have (noticing that
$2\,\hat{\mathfrak h}={\mathfrak h}$)
\begin{eqnarray*}
\psi(\vec e)&=&\sum_{i=1}^{n}\psi(\vec \varepsilon_i) + \psi(\vec e_n) \mod 2 \pi i \\ %
&=& \sum_{i=1}^{n} \Big( \hat{\mathfrak h}(\phi(Y_{i-1})) + \hat{\mathfrak h}(\phi(Y_{i})) + %
\textstyle{\sum_{X\in\Omega^{-}(\vec \varepsilon_i)}}\,{\mathfrak h}(\phi(X)) \Big) + \psi(\vec e_n) \mod 2 \pi i \\ %
&=& \hat{\mathfrak h}(\phi(Y_0)) + \hat{\mathfrak h}(\phi(Y_n)) + \textstyle{\sum_{i=1}^{n-1}}\,{\mathfrak h}(\phi(Y_i)) + \\ %
& & +\sum_{i=1}^{n}\sum_{X\in\Omega^{-}(\vec \varepsilon_i)}{\mathfrak h}(\phi(X)) + \psi(\vec e_n) \mod 2 \pi i. %
\end{eqnarray*}
Then (\ref{eqn:partial sum 2}) follows by letting $n \rightarrow \infty$: %
\begin{eqnarray*}
\psi(\vec e) &=&\hat{\mathfrak h}(\phi(Y_0))+\lim_{n \rightarrow
\infty}\hat{\mathfrak
h}(\phi(Y_n))+\textstyle{\sum_{i=1}^{\infty}}\,{\mathfrak
h}(\phi(Y_i))+\\
& &+\sum_{i=1}^{\infty}\sum_{X\in\Omega^{-}(\vec
\varepsilon_i)}{\mathfrak h}(\phi(X))+\lim_{n \rightarrow
\infty}\psi(\vec e_n) \mod 2 \pi i \\
&=&\hat{\mathfrak
h}(\phi(Y_0))+\textstyle{\sum_{X\in\Omega^{-}(\vec
e)}}\,{\mathfrak h}(\phi(X))+\hat{\mathfrak h}(\phi(X_0)) \mod 2
\pi
i \\
&=&\sum_{X \in \Omega^{0}(e)} \hat {\mathfrak h}(\phi(X))+\sum_{X
\in \Omega^{-}(\vec e)}{\mathfrak h}(\phi(X)) \mod 2 \pi i
\end{eqnarray*}
since it can be easily checked that
\begin{eqnarray}
\lim_{n \rightarrow \infty}\hat{\mathfrak h}(\phi(Y_n))=0,
\end{eqnarray}
\begin{eqnarray}
\lim_{n \rightarrow \infty}\psi(\vec e_n)=\hat{\mathfrak
h}(2)=\hat{\mathfrak h}(\phi(X_0)).
\end{eqnarray}

Proposition \ref{lem:convergence around pm 2} is thus proved.
\square


\vskip 10pt

{\bf Proof of Theorems {\rm\ref{thm:AA}} and {\rm\ref{thm:BB}}.}
\, The\, proofs\, for\, Theorems \ref{thm:AA} and \ref{thm:BB}\,
now pretty much follow in the same way as the proof for the cases
satisfying the BQ-conditions. We construct the finite tree $T$ as
in Step 2; the only difference here is that the circular set $C(T)
\subset \vec E$ may have some directed edges $\vec e$ such that
$\log^{+}|\phi|$ does not have lower Fibonacci bounds on
$\Omega^{-}(\vec e)$. These edges are adjacent to regions on which
$\phi$ takes the value $\pm 2$. We then apply Lemma
\ref{lem:convergence around pm 2} to complete the proofs. \square

\vskip 10pt

\begin{lem}\label{lem:x0=2cos} Suppose a nontrivial $\mu$-Markoff map
$\phi \in {\bf \Phi}_\mu$ is invariant under the action of a
parabolic subgroup $\langle H \rangle$ of $\PSLTwoZ$ where $H \in
\PSLTwoZ$ fixes region $X_0 \in \Omega$ and is conjugated to
$\Big(\begin{array}{cc} 1 & p \\ 0 & 1
\end{array}\Big)$ with $p>1$. Then
$\phi(X_0)=2\cos\frac{q\pi}{p}$ for an integer $q$ with $0 \le q < p$. %
\end{lem}

\begin{pf} Let $Y_n, n \in \mathbf Z$ be the bi-infinite sequence of
neighboring regions of $X_0$. Let
$x_0=\phi(X_0)=\lambda+\lambda^{-1}$. If $x_0 \neq \pm 2$, then
there exist $A,B\in \mathbf C$ such that $AB=(x_0^2-\mu)/(x_0^2-4)$ and %
$y_n=\phi(Y_n)=A\lambda^n+B\lambda^{-n}$. Now we may suppose $H$
acts on $\Omega$ by $H(Y_n)=Y_{n+p}$. Hence we have
$y_{-p}=y_0=y_p$, that is,
\begin{eqnarray}\label{eqn:=A+B=}
A\lambda^{-p}+B\lambda^{p}=A+B=A\lambda^{p}+B\lambda^{-p}.
\end{eqnarray}
Since $\phi$ is nontrivial, $A$ and $B$ are not both $0$. Then it
is easy to derive from (\ref{eqn:=A+B=}) that we always have
$\lambda^{n}=1$ or $\lambda^{n}=-1$, and the latter case only
occurs when $A+B=0$. Hence $\phi(X_0)$ has the desired form.
\end{pf}

\begin{rmk} From the above proof we know that if $y_n \neq 0$ for all $n \in
\mathbf Z$, then we will have $x_0=\phi(X_0)=2\cos\frac{q\pi}{p}$
with $0 \le q < p$.
\end{rmk}

\vskip 10pt

{\bf Proof of Theorem {\rm \ref{thm:CC}}.}\,\, By Lemma
\ref{lem:x0=2cos}, in this case $\phi(X_0)=2\cos (2q\pi/p)$ for
some integers $p > 1$ and $0 < q < p$ with ${\rm gcd}(p,q)=1$. If
$\phi(X_0)\neq \pm \sqrt{\mu}$, we can follow the proof for the
case where $\phi$ satisfies the BQ-conditions and use Lemma
\ref{lem:convergence around pm 2} to deal with the directed edges
$\vec e$ in the circular set for which $\log^{+}|\phi|$ does not
have lower Fibonacci bounds on $\Omega^{-}(\vec e)$. The remaining
case to take care of is when $\phi(X_0)=\pm \sqrt \mu$ (in this
case we cannot fully apply Theorem \ref{thm:partialsum} since
$\Psi_{\mu}(x,y,z)$ is not defined when $x^2$ or $y^2$ is equal to
$\mu$). This is dealt in the following lemma.

\vskip 5pt

\begin{lem}\label{lem:x^2=mu} Suppose $(\phi(X_0))^2=\mu$ in Theorem
{\rm\ref{thm:CC}}. Let $(Y_n), n \in {\mathbf Z}$ be the
bi-infinite sequence of neighboring regions of $X_0$. Let $\vec
\varepsilon_n, n \in {\mathbf Z}$ be the directed edge on
$\varepsilon_n=Y_{n-1} \cap Y_n$ with its head on $X_0$. Then
\begin{eqnarray}\label{eqn:sum psi =0}
\textstyle{\sum_{k=1}^{2p}}\,\psi(\vec \varepsilon_{k})=0 \mod
2\pi i,
\end{eqnarray}
where $p>1$ is the index of $\langle H \rangle$ as a subgroup in a
maximal parabolic subgroup of $\PSLTwoZ$, that is, $H$ is
conjugated to $\Big(\begin{array}{cc}1&p\\0&1\end{array}\Big)$.
\end{lem}

\begin{pf} In this case $\phi(X_0)=2\cos(q\pi/p)$ and
$\mu=4\cos^{2}(q\pi/p)$. It can be shown that
\begin{eqnarray}
\textstyle{\nu = \cosh^{-1}\big(1-{\mu}/{2}\big) = \big(\pi -
{2q\pi}/{p}\big) i \quad {\rm or}\quad \big({2q\pi}/{p}-\pi\big)i,} %
\end{eqnarray}
depending on whether $q \le p/2$ or $q \ge p/2$. Without loss of
generality, we may assume $q \le p/2$. Hence\, $\nu= \big(\pi -
\frac{2q\pi}{p}\big) i$\, and\, $-e^\nu =e^{-i\frac{2q\pi}{p}}$.

\vskip 5pt

Following our convention, let $y_n=\phi(Y_n), n \in {\mathbf Z}$.
By (\ref{eqn:x^2+y^2+z^2-xyz=mu}) we have either
\begin{eqnarray}\label{eqn:exp i pi over p}
\textstyle{\frac{y_j}{y_{j-1}}}=e^{i\frac{q\pi}{p}} \,\,\,\,\,\,
\textrm{ for all } j \in {\mathbf Z}
\end{eqnarray}
or
\begin{eqnarray}\label{eqn:exp -i pi over p}
\textstyle{\frac{y_j}{y_{j-1}}}=e^{-i\frac{q\pi}{p}} \,\,\,\,\,\,
\textrm{ for all } j \in {\mathbf Z}.
\end{eqnarray}
In particular, \begin{eqnarray}y_{2p+j}=y_{j} \,\,\,\,\,\,
\textrm{ for all } j \in {\mathbf Z}. \end{eqnarray}

\vskip 5pt

Recall that
\begin{eqnarray}
\psi(\vec \varepsilon_{j})=\log
\frac{1+(e^{\nu}-1)(x_0/y_{j-1}y_{j})}{\sqrt{1-\mu/y_{j-1}^2}\sqrt{1-\mu/y_{j}^2}}.
\end{eqnarray}
Hence the desired equality (\ref{eqn:sum psi =0}) holds if and
only if
\begin{eqnarray}\label{eqn:prod=1}
\prod_{j=1}^{2p}\frac{1+(e^{\nu}-1)(x_0/y_{j-1}y_{j})}{\sqrt{1-\mu/y_{j-1}^2}\sqrt{1-\mu/y_{j}^2}}=1;
\end{eqnarray}
if and only if
\begin{eqnarray}\label{eqn:varprod=1}
\prod_{j=1}^{2p}\frac{y_{j-1}y_{j}+(e^{\nu}-1)x_0}{y_{j}^2-\mu}=1;
\end{eqnarray}
if and only if 
\begin{eqnarray}\label{eqn:redprod=1}
\frac{y_{0}y_{1}+(e^{\nu}-1)x_0}{y_{1}y_{2}+(e^{-\nu}-1)x_0} \,
\frac{y_{2}y_{3}+(e^{\nu}-1)x_0}{y_{3}y_{4}+(e^{-\nu}-1)x_0} \,
\cdots \,
\frac{y_{2p-2}y_{2p-1}+(e^{\nu}-1)x_0}{y_{2p-1}y_{2p}+(e^{-\nu}-1)x_0}
=1
\end{eqnarray}
since it can be easily checked that
\begin{eqnarray}
[\,y_{j-1}y_{j}+(e^{\nu}-1)x_0]\,[\,y_{j-1}y_{j}+(e^{-\nu}-1)x_0]=(y_{j-1}^2-\mu)(y_{j}^2-\mu).
\end{eqnarray}

\vskip 5pt

{\bf Case 1.}\,\, Suppose (\ref{eqn:exp i pi over p}) holds. Then
$e^{-i\frac{q\pi}{p}} y_{j}=y_{j-1}$ for all $j \in {\mathbf Z}$.
Hence
\begin{eqnarray}
(-e^{\nu})y_{j}y_{j+1}=y_{j-1}y_{j}\,\,\,\,\,\,
\textrm{ for all } j \in {\mathbf Z},
\end{eqnarray}
and we have
\begin{eqnarray*}
\textrm{LHS of} \,
(\ref{eqn:redprod=1})&=&\frac{y_{0}y_{1}+(e^{\nu}-1)x_0}{(-e^{\nu})(y_{1}y_{2}+(e^{-\nu}-1)x_0)}
\,
\frac{y_{2}y_{3}+(e^{\nu}-1)x_0}{(-e^{\nu})(y_{3}y_{4}+(e^{-\nu}-1)x_0)}
\, \\ & & \cdots \,
\frac{y_{2p-2}y_{2p-1}+(e^{\nu}-1)x_0}{(-e^{\nu})(y_{2p-1}y_{2p}+(e^{-\nu}-1)x_0)}\\
&=&1.
\end{eqnarray*}

\vskip 5pt

{\bf Case 2.}\,\, Suppose (\ref{eqn:exp -i pi over p}) holds. Then
$e^{-i\frac{q\pi}{p}} y_{j-1}=y_{j}$ for all $j \in {\mathbf Z}$.
Hence
\begin{eqnarray}
(-e^{\nu})y_{j-1}y_{j}=y_{j}y_{j+1}\,\,\,\,\,\, \textrm{ for all }
j \in {\mathbf Z},
\end{eqnarray}
and we have (noticing that $y_{2p+j}=y_{j}$)
\begin{eqnarray*}
\textrm{LHS of} \,
(\ref{eqn:redprod=1})&=&\frac{y_{0}y_{1}+(e^{\nu}-1)x_0}{(-e^{\nu})[\,y_{1}y_{2}+(e^{-\nu}-1)x_0]}
\,
\frac{y_{2}y_{3}+(e^{\nu}-1)x_0}{(-e^{\nu})[\,y_{3}y_{4}+(e^{-\nu}-1)x_0]}
\,\\ & & \cdots \,
\frac{y_{2p-2}y_{2p-1}+(e^{\nu}-1)x_0}{(-e^{\nu})[\,y_{2p-1}y_{2p}+(e^{-\nu}-1)x_0]}\\
&=& \frac{y_{0}y_{1}+(e^{\nu}-1)x_0}{y_{2}y_{3}+(e^{\nu}-1)x_0} \,
\frac{y_{2}y_{3}+(e^{\nu}-1)x_0}{y_{4}y_{5}+(e^{\nu}-1)x_0} \,
\frac{y_{4}y_{5}+(e^{\nu}-1)x_0}{y_{6}y_{7}+(e^{\nu}-1)x_0}
\\ & & \cdots \,
\frac{y_{2p-4}y_{2p-3}+(e^{\nu}-1)x_0}{y_{2p-2}y_{2p-1}+(e^{\nu}-1)x_0}
\,
\frac{y_{2p-2}y_{2p-1}+(e^{\nu}-1)x_0}{y_{2p}y_{2p+1}+(e^{\nu}-1)x_0}\\
&=&1.
\end{eqnarray*}

This proves Lemma \ref{lem:x^2=mu} and hence Theorem \ref{thm:CC}.
\end{pf}

{\bf Proof of Theorem {\rm \ref{thm:D}}.}\,\, Let us write
$X_{\infty}$ for $X_{\frac10}$, and $X_n$ for $X_{\frac{n}{1}}, n
\in {\mathbf Z}$. (The reader is referred to Figure
\ref{fig:Farey} for an identification of $\Omega$ with ${\mathbf
Q} \cap \{\frac10\}$.) By our convention, we write
$x_{\infty}=\phi(X_{\infty})$ and $x_n=\phi(X_n)$ for $n \in
{\mathbf Z}$. Since $\Omega_\phi(2)=\{X_{\infty}\}$, we know that
$|x_n| > 2$ for all $n \in {\mathbf Z}$.

\vskip 3pt

First, we note that $x_n, n \in {\mathbf Z}$ are bounded.
Actually, let $x_{\infty}=\lambda + {\lambda}^{-1}$; then
$|\lambda|=1$ (since $x_{\infty} \in (-2,2)$), and there exist
$A,B \in \mathbf C$ such that
$AB=(x_{\infty}^2-\mu)/(x_{\infty}^2-4)$ and
$x_n=A\lambda^n+B\lambda^{-n}$ for all $n \in {\mathbf Z}$. Hence
$|x_n| \le |A|+|B|$ for all $n \in {\mathbf Z}$.

\vskip 3pt

Next, we show that $x_n^2, n \in {\mathbf Z}$ are bounded away
from $\mu$. Actually, since
$$x_{n-1}^2+x_n^2+x_{\infty}^2-x_{n-1}x_n x_{\infty}=\mu \quad
{\rm and} \quad x_{n-1}+x_{n+1}=x_n x_{\infty},$$ we have
$$x_n^2-\mu=x_{n-1}(x_n x_{\infty}-x_{n-1})-x_{\infty}^2=x_{n-1}x_{n+1}-x_{\infty}^2,$$
hence
$$|x_n^2-\mu| \ge |x_{n-1}|\,|x_{n+1}|-x_{\infty}^2 > 4-x_{\infty}^2 > 0.$$

\vskip 3pt

{\bf Case I.}\,\, Suppose $x_{\infty}^2 \neq \mu$.

\vskip 3pt

Let $e_n, n \in {\mathbf Z}$ be the edge $X_{\infty} \cap X_n$,
and let $\vec e_n$ be the directed edge on $e_n$, directed from
$e_{n-1}$ to $e_{n+1}$. Then by (\ref{eqn:psi(e)})
$$\psi(\vec e_n)=\log\frac{x_{\infty}x_n+(e^\nu-1)x_{n+1}}{\sqrt{x_{\infty}^2-\mu}\sqrt{x_n^2-\mu}}.$$

Now we claim that $\psi(\vec e_n), n \in {\mathbf Z}$ are bounded.
Actually, it is clear that both $x_{\infty}x_n+(e^\nu-1)x_{n+1}$
and $\sqrt{x_{\infty}^2-\mu}\sqrt{x_n^2-\mu}$ are bounded;
furthermore, the latter is bounded away from $0$ as just shown.
That $x_{\infty}x_n+(e^\nu-1)x_{n+1}$ is bounded away from $0$
follows from the following identity:
\begin{eqnarray}
[\,x_{\infty}x_n+(e^\nu-1)x_{n+1}][\,x_{\infty}x_n+(e^{-\nu}-1)x_{n+1}]=(x_{\infty}^2-\mu)(x_n^2-\mu).
\end{eqnarray}
This proves the claim.

\vskip 3pt

We further claim that $\hat{\mathfrak h}(x_n), n \in {\mathbf Z}$
are also bounded. Recall that by definition
\begin{eqnarray*}
\hat{\mathfrak
h}(x_n)=\log\frac{1+(e^\nu-1)h(x_n)}{\sqrt{1-\mu/x_n^2}}.
\end{eqnarray*}
It is clear that $1+(e^\nu-1)h(x_n), n \in {\mathbf Z}$ are
bounded. That they are bounded away from $0$ follows from the
following identity:
\begin{eqnarray}
[1+(e^\nu-1)h(x_n)][1+(e^{-\nu}-1)h(x_n)]=1-\mu/x_n^2,
\end{eqnarray}
together with the fact shown earlier that $x_n^2,  n \in {\mathbf
Z}$ are bounded away from $\mu$. This proves the claim.

\vskip 3pt

Let $\varepsilon_n, n \in {\mathbf Z}$ be the edge $X_{n-1} \cap
X_n$ and $\vec \varepsilon_n, n \in {\mathbf Z}$ be the directed
edge $(X_{n-1},X_n;X_{\frac{2n-1}{2}} \rightarrow X_\infty)$ on
$\varepsilon_n$.

\vskip 3pt

By Proposition \ref{thm:partialsum} we have for each $n \in
{\mathbf Z}$,
\begin{eqnarray}
\psi(\vec \varepsilon_n)=\hat{\mathfrak h}(x_{n-1})+\hat{\mathfrak
h}(x_{n})+\textstyle{\sum_{r\in{\mathbf Q}\cap(n-1,n)}}{\mathfrak
h}(x_{r}) \mod 2 \pi i,
\end{eqnarray}
and hence (noticing that $2\hat{\mathfrak h}={\mathfrak h} \mod 2
\pi i$)
\begin{eqnarray}\label{eqn:sum psi vare 1}
\textstyle{\sum_{n=1}^{N}}\psi(\vec \varepsilon_n)=\hat{\mathfrak
h}(x_{0})+\hat{\mathfrak h}(x_{N})+\textstyle{\sum_{r\in{\mathbf
Q}\cap(0,N)}}{\mathfrak h}(x_{r}) \mod 2 \pi i.
\end{eqnarray}

Note that for $N \ge 1$, $C_N=\{\vec e_0, \vec \varepsilon_1,
\cdots, \vec \varepsilon_N, -\vec e_N\} \subset \vec E(\Sigma)$ is
a circular set. Actually, $C_N=C(T_N)$ for the subtree
$T_N=\bigcup_{k=1}^{N-1}e_k$ of $\Sigma$ (when $N=1$, $T$ is the
vertex $v=X_0 \cap X_1 \cap X_\infty$). Hence by (\ref{eqn:sum psi
ii}) and (\ref{eqn:sum psi i}), we have
\begin{eqnarray}\label{eqn:sum psi vare 2}
\textstyle{\sum_{n=1}^{N}}\psi(\vec \varepsilon_n)=\nu-\psi(-\vec
e_N)-\psi(\vec e_0)=\psi(\vec e_N)-\psi(\vec e_0) \mod 2 \pi i.
\end{eqnarray}
It follows from (\ref{eqn:sum psi vare 1}) and (\ref{eqn:sum psi
vare 2}) that
\begin{eqnarray}\label{eqn:sum frakh (0,N]}
\textstyle{\sum_{r\in{\mathbf Q}\cap(0,N]}}\,{\mathfrak h}(x_{r}) %
=\hat{\mathfrak h}(x_{N})-\hat{\mathfrak h}(x_{0})+ \psi(\vec e_{N})-\psi(\vec e_0) \mod 2 \pi i. %
\end{eqnarray}
Since it is shown that $\psi(\vec e_N)$ and $\hat{\mathfrak
h}(x_N)$, $N \ge 1$ are all bounded, we have
\begin{eqnarray}\label{eqn:limRe(sum/N)=}
\lim_{N\rightarrow\infty}\Re
\bigg(\frac{\textstyle{\sum_{r\in{\mathbf
Q}\cap(0,N]}}\,{\mathfrak h}(x_{r})}{N}\bigg)=0.
\end{eqnarray}
\noindent This proves (\ref{eqn:phi(X0)in(-2,2)averaging}) for the
case where $x_{\infty}^2 \neq \mu$.

\vskip 5pt

{\bf Case II.}\,\, Suppose $x_{\infty}^2 = \mu$.

\vskip 3pt

In this case $\psi(\vec e_n)(=\infty)$ is not defined. We have to
evaluate $\textstyle{\sum_{r\in{\mathbf Q}\cap(0,N]}}\,{\mathfrak
h}(x_{r})$ $({\rm mod}\, 2 \pi i)$ directly. We shall assume the
notation used in Case I.

\vskip 3pt

Note that $\nu$ is purely imaginary: $\nu = \cosh^{-1}(1-\mu/2)
\in (0, \pi)i$. Since $AB=0$, we may assume that $x_n=A\lambda^n$
without loss of generality (recall that
$x_\infty=\lambda+\lambda^{-1}$). It is easy to check that we have
either $-e^\nu=\lambda^{-2}$ or $-e^\nu=\lambda^{2}$; and hence
correspondingly, either (\ref{eqn:-2}) or (\ref{eqn:+2}) below
holds:
\begin{eqnarray}
&&\hspace{-20pt}\textstyle{x_{n-1}x_n+(e^{\nu}-1)x_\infty=(-e^\nu)[x_n
x_{n+1}+(e^{-\nu}-1)x_\infty]} \quad {\rm for \,\, all \,\,\,} n
\in \mathbf Z, \label{eqn:-2} \\ && \nonumber \\
&&\hspace{-20pt}\textstyle{x_{n}x_{n+1}+(e^{\nu}-1)x_\infty=(-e^\nu)[x_{n-1}
x_n+(e^{-\nu}-1)x_\infty]} \quad {\rm for \,\, all\,\,\,} n \in
\mathbf Z.\label{eqn:+2}
\end{eqnarray}

When $N$ is even, we have, with similar calculations as in the
proof of Lemma \ref{lem:x^2=mu}, that either (when $-e^\nu=\lambda^{-2}$ hence (\ref{eqn:-2}) holds) %
\begin{eqnarray}\label{eqn:sum psi vare case II}
\sum_{n=1}^{N}\psi(\vec \varepsilon_n)&=&\log\frac{\sqrt{1-\mu/x_N^2}}{\sqrt{1-\mu/x_0^2}} + %
\log\frac{\prod_{n=1}^{N}[x_{n-1}x_n+(e^\nu-1)x_\infty]}{\prod_{n=1}^{N}(x_n^2-\mu)}  \mod 2\pi i \nonumber \\
&=&\log\frac{\sqrt{1-\mu/x_N^2}}{\sqrt{1-\mu/x_0^2}} + \log \bigg(
\frac{x_{0}x_{1}+(e^{\nu}-1)x_\infty}{x_{1}x_{2}+(e^{-\nu}-1)x_\infty} %
\frac{x_{2}x_{3}+(e^{\nu}-1)x_\infty}{x_{3}x_{4}+(e^{-\nu}-1)x_\infty} \nonumber \\ %
& &\hspace{65pt} \, \cdots \, \frac{x_{N-2}x_{N-1}+(e^{\nu}-1)x_\infty}{x_{N-1}x_{N}+
   (e^{-\nu}-1)x_\infty} \bigg) \mod 2\pi i \nonumber \\ %
&=&\log\frac{\sqrt{1-\mu/x_N^2}}{\sqrt{1-\mu/x_0^2}} + \frac{N}{2}(\nu+\pi i) \mod 2 \pi i %
\end{eqnarray}
or (when $-e^\nu=\lambda^{2}$ hence (\ref{eqn:+2}) holds) %
\begin{eqnarray}\label{eqn:sum psi vare case II or}
\hspace{-25pt}\sum_{n=1}^{N}\psi(\vec \varepsilon_n)&=&
\log\frac{\sqrt{1-\mu/x_N^2}}{\sqrt{1-\mu/x_0^2}}+\frac{N}{2}(\nu+\pi i) \nonumber \\ %
& & \hspace{72pt} + \,\log\frac{x_{0}x_{1}+(e^{\nu}-1)x_\infty}{x_{N}x_{N+1}+ (e^{\nu}-1)x_\infty} \mod 2 \pi i. %
\end{eqnarray}
When $N$ is odd, we have similar expressions with one more term in each case since %
$$\textstyle{\sum_{n=1}^{N}}\psi(\vec \varepsilon_n)=\psi(\vec \varepsilon_N)+\textstyle{\sum_{n=1}^{N-1}}
\psi(\vec \varepsilon_n)$$ and $N-1$ is even. Hence we still have
(\ref{eqn:limRe(sum/N)=}) in this case, noticing that (a) $\nu$ is
purely imaginary, (b) $x_N^2, N \ge 1$ are bounded away from $\mu$
and $0$, and (c) $\psi(\vec \varepsilon_N), N \ge 1$ are bounded.
This proves (\ref{eqn:phi(X0)in(-2,2)averaging}) for the case
where $x_{\infty}^2 = \mu$.

\vskip 5pt

The proof of (\ref{eqn:phi(X0)in(-2,2)asymptotic}) is similar. 
Recall that $x_\infty=\lambda+\lambda^{-1} \in (-2,2)$, hence
$|\lambda|=1$. Thus there exists an increasing sequence $(N_k)_{k
\ge 1}$ of positive integers such that $\lim_{k \rightarrow
\infty}\lambda^{N_k}=1$. Therefore we have $\lim_{k \rightarrow
\infty}x_{N_k}=x_0$, $\lim_{k \rightarrow \infty}\psi(\vec
e_{N_k})=\psi(\vec e_0)$ and $\lim_{k \rightarrow
\infty}\hat{\mathfrak h}(x_{N_k})=\hat{\mathfrak h}(x_0)$. Hence
if $x_\infty^2\neq \mu$, we have from (\ref{eqn:sum frakh (0,N]})
(with $N=N_k$) that
\begin{eqnarray}
\textstyle{\sum_{r\in{\mathbf Q}\cap(0,N_k]}}\,{\mathfrak h}(x_{r}) %
=\hat{\mathfrak h}(x_{N_k})-\hat{\mathfrak h}(x_{0})+ \psi(\vec e_{N_k})-\psi(\vec e_0) \mod 2 \pi i. %
\end{eqnarray}
It follows that $\lim_{k \rightarrow
\infty}\textstyle{\sum_{r\in{\mathbf Q}\cap(0,N_k]}}\,{\mathfrak
h}(x_{r})=0$ (${\rm mod}\, 2 \pi i$)\, in the case $x_\infty^2\neq
\mu$.

\vskip 3pt

Now suppose $x_\infty^2 = \mu$. We may assume that $N_k, k \ge 1$
are all even. Then in the case when $-e^\nu=\lambda^{-2}$, we have
from (\ref{eqn:sum psi vare 1}) and (\ref{eqn:sum psi vare case
II}) (with $N=N_k$) that
\begin{eqnarray}
\textstyle{\sum_{r\in{\mathbf Q}\cap(0,N_k]}}\,{\mathfrak h}(x_{r}) %
&=& \hat{\mathfrak h}(x_{N_k})-\hat{\mathfrak h}(x_{0})
+\log\Big[\sqrt{1-\mu/x_{N_k}^2}\, \Big / \sqrt{1-\mu/x_0^2}\,\,\Big] \nonumber \\ %
& & +\, \log\big[(-e^\nu)^{N_k/2}\big]  \mod 2 \pi i. %
\end{eqnarray}
We have again $\lim_{k \rightarrow
\infty}\textstyle{\sum_{r\in{\mathbf Q}\cap(0,N_k]}}\,{\mathfrak
h}(x_{r})=0$ (${\rm mod}\, 2 \pi i$)\, in this case, noticing that
$$\lim_{k \rightarrow
\infty}\log\big[(-e^\nu)^{N_k/2}\big]= \lim_{k \rightarrow
\infty}\big[(\lambda^{-2})^{N_k/2}\big]=\log 1 = 0 \mod 2 \pi i.$$ %
The case when $-e^\nu=\lambda^{2}$ can be proved similarly, this
time using (\ref{eqn:sum psi vare 1}) and (\ref{eqn:sum psi vare
case II or}). This finishes the proof of
(\ref{eqn:phi(X0)in(-2,2)asymptotic}) and hence also the proof of
Theorem \ref{thm:D}. \square


 \vskip 18pt
\section{\bf Geometric interpretations and examples}\label{sec:geometric}
 \vskip 10pt

In this section we give some examples of (generalized) Markoff
maps $\phi \in \Phi_{\mu}$ which are extremal with respect to the
conditions, namely where $\phi(X)=\pm 2$ for a finite number of $X
\in \Omega$ (or $[X] \in \Omega/\langle H\rangle$ in the case of
Theorem \ref{thm:CC}), and the geometric situations when they
occur.

\begin{example}\label{ex:AA}
We start with examples for Theorem \ref{thm:AA}. The Maskit
embedding of ${\mathcal T}_{1,1}$ parameterized by the complex
number $\zeta$, with corresponding Kleinian group $\Gamma(\zeta)$
as studied by Linda Keen and Caroline Series in
\cite{keen-series1993t} corresponds to the $0$-Markoff triple $(2,
-i\zeta, -i(\zeta + 2))$ and the corresponding 0-Markoff map
$\phi_{\zeta}\in \Phi_0$ (they used the parameter $\mu$ instead of
$\zeta$; here we use $\zeta$ to avoid confusion as $\mu$ is
already used for the generalized Markoff maps). These satisfy the
conditions of Theorem \ref{thm:AA} since the number of simple
closed geodesics on the hyperbolic 3-manifold ${\bf
H}^3/\Gamma(\zeta)$ of length $\le K$, for any $K>0$ is finite.
For example, if we take $\zeta$ such that its imaginary part
${\Im}\zeta \ge 2$, then it is easy to see that the conditions
will be satisfied. These examples correspond to points on the
boundary of quasifuchsian space ${\mathcal QF}_{1,1}$ for which
McShane's identity still holds. More generally, the same is true
even when we look at the cusp points on the boundary of the Maskit
embedding, these correspond to the case where exactly two curves
are pinched, one on each end of the 3-manifold. On the other hand,
Minsky's results from \cite{minsky1999am}, see also the references
contained therein, imply that for any other boundary point of
${\mathcal QF}_{1,1}$, there is one or two degenerate ends, in
which case there are infinitely many homotopy classes of simple
closed curves on the torus of bounded length, and hence McShane's
identity does not hold. Hence, at least on $\partial ({\mathcal
QF}_{1,1})$, one can determine completely the points where
McShane's identity holds. Now if Bowditch's conjecture is true,
that is, if the set of Markoff maps satisfying the BQ-conditions
are exactly those which come from the quasifuchsian
representations, this would give a complete geometric description
of all representations for which McShane's identity holds, in the
type-preserving case.
\end{example}

\begin{example}\label{ex:BB}
We next give an example of an extremal case for Theorem
\ref{thm:BB}, in the case where $\mu=20$. This is the Markoff map
$\phi \in \Phi_{20}$ corresponding to the Markoff triple
$(-2,-2,-2)$, which arises from the representation coming from the
complete, finite volume hyperbolic structure on the
thrice-punctured sphere. Note that in this case $|\Omega(2)|=3$.
It seems likely that there is a fixed constant $N$ such that
$|\Omega_\phi(2)|\le N$ for all generalized Markoff maps $\phi$
for which McShane's identity holds.

\end{example}

\begin{example}\label{ex:CC}
We finally give some examples for Theorem \ref{thm:CC}. Consider
the Markoff map $\phi \in \Phi_0$ corresponding to the triple
$(0,2,2i)$. It is not difficult to verify that this satisfies the
conditions of Theorem \ref{thm:CC} where the stabilizer $H$ of
$\phi$ is conjugate to $\Big(%
\begin{array}{cc}
  1 & 4 \\
  0 & 1 \\
\end{array}%
\Big) \in \PSLTwoZ$. More generally, so does the map $\phi\in
\Phi_0$ corresponding to the triple $(0, \zeta, \zeta i)$ for
$|\zeta| > 2$. It is not difficult to construct other examples
with other stabilizers, by replacing $0$ with $2 \cos (\pi/n)$,
and the other values appropriately, for example, the $\mu$-Markoff
map $\phi \in \Phi_{\mu}$, where $\mu =4 \cos ^2 (\pi/n)$,
corresponding to the triple $(2 \cos (\pi/n), \zeta,
e^{i\pi/n}\zeta )$ where again $|\zeta| > 2$. It seems that most
of these examples should correspond to representations arising
from hyperbolic three orbifolds. However, we are not sure about
the detailed geometric interpretations of these generalized
Markoff maps.
\end{example}

\vskip 30pt


{}

\end{document}